%% file: paper.tex
\documentclass{lmcs}


\usepackage{hyperref}
\usepackage{amsmath, amsfonts, amsthm}
\usepackage{tikz}
\usepackage{tikz-cd}
\usetikzlibrary{arrows,matrix,decorations.pathmorphing,
  decorations.markings, calc, backgrounds}
\usepackage{mathpartir}
\usepackage{quiver}
\usepackage{xspace}
\usepackage{xcolor}
\usepackage{bbm}
\usepackage{breakurl}
\hypersetup{pdfborder={0 0 100}}

\input{macros}

\begin{document}

\title[Formalising and Computing $\pi_4(\mathbb{S}^3)$ in Cubical
  Agda]{Formalising and Computing the Fourth Homotopy Group of the
  $3$-Sphere in Cubical Agda}

\titlecomment{This paper is an extended version of ``Formalizing
  $\pi_4(\mathbb{S}^3) \cong \mathbb{Z} / 2 \mathbb{Z}$ and Computing
  a Brunerie Number in Cubical Agda'' published in the
  post-proceedings of Logic in Computer Science 2023
  \cite{LICS23}. Some details about what has been added can be found
  in the \textbf{Outline} paragraph below.}

\thanks{This paper is based upon research supported by the Swedish
  Research Council (SRC, Vetenskapsrådet) under Grant
  No.~2019-04545. The research has also received funding from the Knut
  and Alice Wallenberg Foundation through the Foundation's program for
  mathematics.}

\author[A.~Ljungström]{Axel Ljungström\lmcsorcid{0000-0001-6946-0775}}[a]
\author[A.~Mörtberg]{Anders Mörtberg\lmcsorcid{0000-0001-9558-6080}}[a]
\address{Department of Mathematics, Stockholm University, Stockholm, Sweden}
\email{axel.ljungstrom@math.su.se, anders.mortberg@math.su.se}

\begin{abstract}
  \noindent Brunerie's 2016 PhD thesis contains the first synthetic
  proof in Homotopy Type Theory (HoTT) of the classical result that
  the fourth homotopy group of the 3-sphere is
  $\mathbb{Z}/2\mathbb{Z}$. The proof is one of the most impressive
  pieces of synthetic homotopy theory to date and uses a lot of
  advanced classical algebraic topology rephrased
  synthetically. Furthermore, the proof is fully constructive and the
  main result can be reduced to the question of whether a particular
  ``Brunerie number'' $\beta$ can be normalised to $\pm 2$. The
  question of whether Brunerie’s proof could be formalised in a proof
  assistant, either by computing this number or by formalising the
  pen-and-paper proof, has since remained open. In this paper, we
  present a complete formalisation in \CubicalAgda. We do this by
  modifying Brunerie’s proof so that a key technical result, whose
  proof Brunerie only sketched in his thesis, can be avoided. We also
  present a formalisation of a new and much simpler proof that $\beta$
  is $\pm 2$. This formalisation provides us with a sequence of
  simpler Brunerie numbers, one of which normalises very quickly to
  $-2$ in \CubicalAgda, resulting in a fully formalised
  computer-assisted proof that $\pi_4(\mathbb{S}^3) \cong
  \mathbb{Z}/2\mathbb{Z}$.
\end{abstract}

\maketitle

\input{sections/intro}
\input{sections/hott}

\input{sections/homotopygroupsofspheres}
\input{sections/brunerienumber}
\input{sections/brunerieproof}
\input{sections/computation}
\input{sections/conclusion}

\section*{Acknowledgement}
\noindent First and foremost we would like to thank Guillaume Brunerie
for his excellent thesis, the conjecture about the computability of
$\beta$, and for the many discussions about this over the years. We
would also like to thank Thierry Coquand and Simon Huber for the first
attempt to compute the number together with Guillaume using
\systemname{cubical} in December 2014 and everyone else who has tried
to compute the number and contributed ideas to this since. The
\CubicalAgda formalisation relies on many contributions to
\agdaCubical by more people than we can mention, but we are especially
grateful to Evan Cavallo for the Freudenthal suspension theorem as
well as many cool cubical tricks and to Rongji Kang for the
Blakers-Massey theorem.

\bibliographystyle{alphaurl}
\bibliography{refs}

\end{document}

%% file: macros.tex
\newcommand{\teletype}[1]{\ensuremath{\mathtt{#1}}}
\newcommand{\systemname}[1]{\teletype{\color{darkgray}#1}\xspace}

\newcommand{\Agda}{\systemname{Agda}}

\newcommand{\agdaCubical}{\systemname{agda/cubical}}

\newcommand{\CubicalAgda}{\systemname{Cubical} \systemname{Agda}}
\newcommand{\cubicaltt}{\systemname{cubicaltt}}
\newcommand{\cooltt}{\texttt{\textcolor[rgb]{.012,.27,.46}{cool}tt}}

\definecolor{Revolutionary}{RGB}{232,70,68}
\newcommand{\redtt}{\textbf{\teletype{{\color{Revolutionary}red}tt}}}

\usepackage{agda/latex/agda}
\newcommand{\anum}[1]{\AgdaNumber{#1}}

\newcommand{\field}[1]{{\AgdaField{\ensuremath{\mathsf{#1}}}}}

\newcommand{\func}[1]{\AgdaFunction{#1}}

\newcommand{\con}[1]{{\AgdaInductiveConstructor{\ensuremath{\mathsf{#1}}}}}
\usepackage{catchfilebetweentags}
\usepackage{newunicodechar}
\newunicodechar{≃}{\ensuremath{\simeq}}
\newunicodechar{≅}{\ensuremath{\cong}}
\newunicodechar{∈}{\ensuremath{\in}}
\newunicodechar{⁻}{\ensuremath{^{-}}}
\newunicodechar{ᴹ}{\ensuremath{_{M}}}
\newunicodechar{ᴺ}{\ensuremath{_{N}}}
\newunicodechar{≡}{\ensuremath{\equiv}}
\newunicodechar{λ}{\ensuremath{\lambda}}
\newunicodechar{⊎}{\ensuremath{\uplus}}
\newunicodechar{∷}{\ensuremath{::}}
\newunicodechar{ℓ}{\ensuremath{\ell}}
\newunicodechar{ᵢ}{\ensuremath{_i}}
\newunicodechar{⟨}{\ensuremath{\langle}}
\newunicodechar{⟩}{\ensuremath{\rangle}}
\newunicodechar{α}{\ensuremath{\alpha}}
\newunicodechar{β}{\ensuremath{\beta}}
\newunicodechar{θ}{\ensuremath{\theta}}
\newunicodechar{φ}{\ensuremath{\varphi}}
\newunicodechar{ψ}{\ensuremath{\psi}}
\newunicodechar{η}{\ensuremath{\eta}}
\newunicodechar{ε}{\ensuremath{\varepsilon}}
\newunicodechar{ι}{\ensuremath{\iota}}
\newunicodechar{Σ}{\ensuremath{\Sigma}}
\newunicodechar{σ}{\ensuremath{\sigma}}
\newunicodechar{∀}{\ensuremath{\forall}}
\newunicodechar{ℕ}{\ensuremath{\mathbb{N}}}
\newunicodechar{→}{\ensuremath{\to}}
\newunicodechar{⊎}{\ensuremath{\uplus}}
\newunicodechar{⋆}{\ensuremath{*}}
\newunicodechar{¬}{\ensuremath{\lnot}}
\newunicodechar{Θ}{\ensuremath{\theta}}
\newunicodechar{ᴸ}{\ensuremath{{}^{\color{black}\leftarrow}}}  
\newunicodechar{ᴿ}{\ensuremath{{}^{\color{black}\rightarrow}}} 
\newunicodechar{∧}{\ensuremath{\wedge}}
\newunicodechar{∨}{\ensuremath{\vee}}
\newunicodechar{∼}{\ensuremath{\sim}}
\newunicodechar{≢}{\ensuremath{\nequiv}}
\newunicodechar{Π}{\ensuremath{\Pi}}
\newunicodechar{ℤ}{\ensuremath{\mathbb{Z}}}
\newunicodechar{∥}{\ensuremath{\parallel}}
\newunicodechar{∣}{\ensuremath{\mid}}
\newunicodechar{ℚ}{\ensuremath{\mathbb{Q}}}
\newunicodechar{₊}{\ensuremath{_+}}
\newunicodechar{∗}{\ensuremath{\ast}}
\newunicodechar{₁}{\ensuremath{_1}}
\newunicodechar{₂}{\ensuremath{_2}}
\newunicodechar{⊥}{\ensuremath{\bot}}
\newunicodechar{∘}{\ensuremath{\circ}}
\newunicodechar{ρ}{\ensuremath{\rho}}
\newunicodechar{↦}{\ensuremath{\mapsto}}
\newunicodechar{₀}{\ensuremath{_0}}
\newunicodechar{∙}{\ensuremath{\boldsymbol{\cdot}}}
\newunicodechar{Ω}{\ensuremath{\Omega}}
\newunicodechar{§}{\ensuremath{\mathbb{S}}}
\newunicodechar{𝟙}{\ensuremath{\mathbb{1}}}
\newunicodechar{×}{\ensuremath{\times}}
\newunicodechar{♥}{\textnormal{\func{Susp}}}
\newunicodechar{₋}{\func{${}_{-}\!$}}
\newunicodechar{✪}{\func{$\sigma\!$}}
\newunicodechar{♢}{$\sphere{2}$}
\newunicodechar{⋆}{\func{\text{\textasteriskcentered}}}
\newunicodechar{®}{${}^{\text{\color{black}{$3$}}}$}
\newunicodechar{ℏ}{,}
\newunicodechar{⌣}{$\smile$}
\newunicodechar{⋈}{\text{\textasteriskcentered}}
\newunicodechar{▬}{${}^2$}
\newunicodechar{𝕒}{\color{black}{\textnormal{$A$}}}
\newunicodechar{𝕓}{\color{black}{\textnormal{$B$}}}
\newunicodechar{𝕔}{$\star_A$}
\newunicodechar{𝕕}{$\star_B$}
\newunicodechar{∇}{$\nabla$}
\newunicodechar{⊓}{\color{black}{\textnormal{$A$}}}
\newunicodechar{𝕝}{\color{black}{\textnormal{$B$}}}
\newunicodechar{𝕗}{$i^{\vee}$}
\newunicodechar{𝟙}{$\mathbbm{1}$}
\newunicodechar{⊗}{$\mathbb{S}^{{\color{black}2}}$}
\newunicodechar{⊕}{$\mathbb{S}^{{\color{black}3}}$}
\newunicodechar{₃}{${}_3$}
\newunicodechar{∐}{\con{\star_\wedge}}
\newunicodechar{Ą}{$\star_A$}
\newunicodechar{Ę}{$\star_B$}
\newunicodechar{⋀}{\func{$\wedge$}}
\newunicodechar{ₗ}{\textnormal{${}_l$}}
\newunicodechar{ᵣ}{\textnormal{${}_r$}}
\newunicodechar{℻}{\textnormal{\func{$\sigma$}}}

\newunicodechar{０}{\textnormal{$\func{F}_{\color{black}{n,m}}$}}
\newunicodechar{⋂}{\textnormal{$\sphere{\color{black}{n}}$}}
\newunicodechar{⋃}{\textnormal{$\sphere{\color{black}{m}}$}}
\newunicodechar{𝕄}{\textnormal{$\sphere{\color{black}{n+m+1}}$}}

\newunicodechar{◭}{{\color{black}{($f$}}\func{$+^{*}$}\,{\color{black}{$g$)}}}
\newunicodechar{◬}{{\color{black}{($f$}}\func{$+^{\func{Susp}}$}\,{\color{black}{$g$)}}}
\newunicodechar{Ｆ}{\ensuremath{\star_f}}
\newunicodechar{Γ}{\ensuremath{\star_g}}
\newunicodechar{∁}{\ensuremath{\star_C}}
\newunicodechar{►}{\func{$\text{cong}_\star$}\,{\color{black}{$f$}}}
\newunicodechar{▻}{\func{$\text{cong}_\star$}\,{\color{black}{$g$}}}
\newunicodechar{𝕖}{\ensuremath{\star_C}}
\newunicodechar{𝕨}{{\color{black}{($f$}}\, $\cdot_w$ \,{\color{black}{$g$)}}}
\newunicodechar{¹}{{\color{black}{\ensuremath{^{1}}}}}

\newunicodechar{□}{\con{\ensuremath{\square}}}
\newunicodechar{α}{\func{\ensuremath{{\mathbb{T}^2}}}} 

\newunicodechar{β}{\func{\ensuremath{{\mathbb{K}^2}}}} 

\newunicodechar{γ}{\func{{$\gamma$}}}
\newunicodechar{◃}{}
\newunicodechar{▹}{}
\newunicodechar{◁}{$\func{\textsf{Susp}}$}
\newunicodechar{△}{$\,\func{\textsf{Susp}}$}
\newunicodechar{ℜ}{\con{\textsf{inr}}}
\newunicodechar{𝕃}{\con{\textsf{inl}}}
\newunicodechar{⋁}{\func{$\vee$}}
\newunicodechar{✶}{$\star$}
\newunicodechar{▰}{\func{$\sigma$}}

\newcommand{\One}{\func{$\mathbbm{1}$}}

\newunicodechar{ₕ}{\func{\ensuremath{_h}}}
\newunicodechar{ϕ}{\func{\ensuremath{\phi}}}
\newunicodechar{ζ}{\anum{3}\;\;\;\;\;\;\;\;\;}

\newcommand{\Type}{\textnormal{\func{Type}}}
\newcommand{\nType}[1]{\textnormal{${#1}$-\func{Type}}}

\newcommand{\tyPath}[2]{{#1}\,\func{≡}\,{#2}}

\newcommand{\tyEquiv}[2]{{#1}\,\,\func{≃}\,\,{#2}}
\newcommand{\tyIso}[2]{{#1}\,\,\func{$\cong$}\,\,{#2}}

\newcommand{\Square}[5]{
\newcommand{\maxVal}{{(29 * (#1)) / 30}}
\newcommand{\minVal}{{(#1) / 30}}
\begin{center}
    \begin{tikzpicture}
    \draw[->] (0 , \minVal) -- (0 , \maxVal);
    \node[anchor = east] at (0 , #1 / 2) {#2} ;

    \draw[->] (\minVal , #1) -- (\maxVal , #1) ;
    \node[anchor = south] at (#1 / 2 , #1) {#3} ;

    \draw[->] (#1 , \minVal) -- (#1 , \maxVal) ;
    \node[anchor = west] at (#1 , #1 / 2) {#4} ;

    \draw[->] (\minVal , 0) -- (\maxVal , 0) ;
    \node[anchor = north] at (#1 / 2 , 0) {#5} ;

    \end{tikzpicture}
\end{center}
\let\maxVal\undefined
\let\minVal\undefined
}

\newlength{\LETTERheight}
\AtBeginDocument{\settoheight{\LETTERheight}{I}}

\newcommand{\bS}{\mathbb{S}}
\newcommand{\ZZ}{\mathbb{Z}}
\newcommand{\NN}{\mathbb{N}}
\newcommand{\bZ}{\func{ℤ}}

\newcommand{\sphere}[1]{\func{§}^{#1}}

\newcommand{\cohom}[2]{\textnormal{\func{H}}^{#1}\!\left({#2}\right)}

\newcommand{\CP}{\func{$\mathbb{C}P^2$}}
\newcommand{\refl}{\func{refl}}
\newcommand{\inr}[1]{\con{inr}\,{#1}}
\newcommand{\inl}[1]{\con{inl}\,{#1}}
\newcommand{\push}[1]{\con{push}\,{#1}}
\newcommand{\trunc}[1]{\con{∣}\,#1\,\con{∣}}
\newcommand{\truncT}[2]{\func{∥}\,#2\,\func{∥}_{#1}}
\newcommand{\Loop}{\con{loop}}
\newcommand{\base}{\con{base}}
\newcommand{\surf}{\con{surf}}
\newcommand{\north}{\con{north}}
\newcommand{\south}{\con{south}}
\newcommand{\Code}{\func{Code}}
\newcommand{\cofib}{{\textnormal{\func{cofib}}}\,}
\newcommand{\fib}{\textnormal{\func{fib}}\,}
\newcommand{\merid}[1]{\con{merid}\;{#1}}

\newcommand{\ap}[2]{\func{cong}\;{#1}\;{#2}}

\newcommand{\cupprodop}{\func{$\smile$}}
\newcommand{\cupprodkop}{\func{$\smile_k$}}
\newcommand{\cupprodk}[2]{{#1}\;\cupprodkop\;{#2}}
\newcommand{\cupprod}[2]{{#1}\;\cupprodop\;{#2}}
\newcommand{\cupprodaltop}{\cupprodop}
\newcommand{\cupprodalt}[2]{{#1}\;\cupprodaltop\;{#2}}
\newcommand{\brunerie}{\func{$\beta$}}
\newcommand{\abs}[1]{\left | {#1} \right |}
\newcommand{\EM}[1]{\textnormal{\func{K}}_{#1}}
\newcommand{\LoopSpace}[1]{\func{$\Omega$}\,{#1}}
\newcommand{\LoopSpaceN}[2]{\func{$\Omega$}^{#1}\,{#2}}
\newcommand{\plusk}{\func{+$_k$}}
\newcommand{\Smash}{\func{$\wedge$}}
\newcommand{\Susp}[1]{\textnormal{\func{Susp}{$\,#1$}}}
\newcommand{\SuspN}[2]{\textnormal{\func{Susp}${}^{#1}$}{\,#2}}
\newcommand{\JTrunc}{\truncT{3}{\js}}
\newcommand{\sym}[1]{{#1}^{\func{$-1$}}}
\newcommand{\JJoin}[2]{{#1}\func{\text{\textasteriskcentered}}\,{#2}}
\newcommand{\JoinS}{\JJoin{\sphere{1}}{\sphere{1}}}
\newcommand{\prehopf}{\textnormal{\func{h}}}
\newcommand{\hopf}{\textnormal{\func{hopf}}}
\newcommand{\HI}[1]{\textnormal{\func{HI}}\,#1}
\newcommand{\Nabla}{\func{$\nabla$}}
\newcommand{\ivee}{\func{$i^{\vee}$}}
\newcommand{\james}[2]{\ensuremath{\textnormal{\func{J}}{}_{#1}{\,#2}}}
\newcommand{\js}{\james{2}{\sphere{2}}}
\newcommand{\encode}{\textnormal{\func{encode}}}
\newcommand{\decode}{\textnormal{\func{decode}}}
\newcommand{\idfun}{\textnormal{\func{id}}}
\newcommand{\liftL}[1]{\func{$\widehat{\phantom{-}}$}{#1}}

\newcommand{\smashelem}[2]{\con{\langle}{{#1}\,\con{,}\,{#2}}\con{\rangle}}
\newcommand{\smashpt}{\con{\star_\wedge}}

\newcommand{\pushr}{\con{\textsf{push}_{\textnormal{\textsf{r}}}}}

\newcommand{\falsebool}{\con{\textsf{false}}}
\newcommand{\truebool}{\con{\textsf{true}}}

\newcommand{\pinch}{\textnormal{\func{pinch}}}
\newcommand{\etafun}[1]{\textnormal{\func{$\eta_{#1}$-fun}}}
\newcommand{\fdotg}{\ensuremath{f \, \func{$\cdot_w$} \,g}}
\newcommand{\congpt}[2]{\ensuremath{\func{$\text{cong}_\star$}\,{#1}\,{#2}}}

\definecolor{dkblue}{rgb}{0,0.1,0.5}
\definecolor{lightblue}{rgb}{0,0.5,0.5}
\definecolor{dkgreen}{rgb}{0,0.6,0}
\definecolor{dkbrown}{rgb}{0.4,0,0}
\definecolor{dkviolet}{rgb}{0.6,0,0.8}

\newtheorem*{challenge}{Challenge}


%% file: sections/intro.tex
\section{Introduction}
\label{sec:intro}

Homotopy theory originated in algebraic topology, but is by now a
central tool in many branches of modern mathematics, such as algebraic
geometry and category theory. One of the central notions of study in
homotopy theory is that of the \emph{homotopy groups} of a space $X$,
denoted $\pi_n(X)$. These groups constitute a topological invariant,
making them a powerful tool for establishing whether two given spaces
can or cannot be homotopy equivalent. The first two such groups of a
space are easy to understand: $\pi_0(X)$ characterises the connected
components of $X$ and $\pi_1(X)$ is the fundamental group, i.e.\ the
group of equivalence classes consisting of the loops contained in $X$
up to homotopy. This idea generalises to higher values of $n$, for
which $\pi_n(X)$ consists of $n$-dimensional loops up to homotopy. For
many spaces, these groups tend to become increasingly esoteric and
difficult to compute for large $n$. This is true also for
seemingly tame spaces like spheres, for which
$\pi_n(\mathbb{S}^m)$ in general is highly irregular when $n > m \geq
2$.\footnote{See~\cite[Figure~2.1]{Brunerie16} for a table of
  $\pi_n(\mathbb{S}^m)$ for small $n$ and $m$.}  This paper concerns
the first computer formalisation of the classical result that
$\pi_4(\mathbb{S}^3) \cong \ZZ / 2 \ZZ$, a result which is
particularly interesting because it gives the whole first stable stem
of homotopy groups of spheres, i.e.\ $\pi_{n+1}(\mathbb{S}^n)$ for $n \geq
3$. The fact that $\pi_4(\mathbb{S}^3) \cong \ZZ / 2 \ZZ$ was proved
already in the 1930's by Pontryagin using
cobordism theory, but we instead follow the synthetic approach to
homotopy theory developed in Homotopy Type Theory (HoTT) and
popularised by the HoTT Book~\cite{HoTT13}. In this new approach to
homotopy theory, spaces are represented directly as (higher inductive)
types and homotopy groups are computed using Voevodsky's univalence
axiom~\cite{Voevodsky10cmu}. This gives a logical approach to homotopy
theory, suitable for computer formalisation in proof assistants based
on type theory, while also making it possible to interpret results
in any suitably structured $(\infty,1)$-topos \cite{Shulman19}.

The basis for our formalisation is the 2016 PhD thesis of Brunerie
\cite{Brunerie16} which contains the first synthetic proof in HoTT
that $\pi_4(\bS^3) \cong \ZZ/2\ZZ$. The proof is one of the most
impressive pieces of synthetic homotopy theory to date and uses
advanced machinery from classical algebraic topology developed
synthetically, including the symmetric monoidal structure of smash
products, (integral) cohomology rings, the Mayer-Vietoris and Gysin
sequences, the Hopf invariant, Whitehead products, etc. The
formalisation of Brunerie's proof has since remained open, primarily
due to the highly technical nature of some of the details. In this
paper, we will present such a formalisation in \CubicalAgda
\cite{VezzosiMortbergAbel19}, a \emph{cubical} extension of the \Agda
proof assistant \cite{Agda} with native support for computational
univalence and higher inductive types (HITs).

In addition to being a very impressive proof in synthetic homotopy
theory, Brunerie's proof is particularly interesting as it is fully
constructive. The proof consists of two parts, with the first one
culminating in Chapter 3 with the definition of a number $\beta : \ZZ$
such that $\pi_4(\bS^3) \cong \ZZ/\beta\ZZ$. Since then, this $\beta$
has been commonly referred to as the \emph{Brunerie number}. Brunerie
writes the following about it:

\begin{quote}
  \emph{This result is quite remarkable in that even though it is a
    constructive proof, it is not at all obvious how to actually
    compute this [$\beta$]. At the time of writing, we still haven’t
    managed to extract its value from its
    definition. \cite[Page~85]{Brunerie16}}\end{quote}

In fact, \cite[Appendix B]{Brunerie16} contains a complete and concise
definition of $\beta$ as the image of $1$ under a sequence of $12$ maps:

  \[
\begin{tikzcd}
  \ZZ \arrow[r,"n\mapsto\mathsf{loop}^n"] &
  \Omega(\mathbb{S}^1) \arrow[r,"\Omega \varphi_{\mathbb{S}^1}"] &
  \Omega^2(\mathbb{S}^2) \arrow[r,"\Omega^2 \varphi_{\mathbb{S}^2}"] &
  \Omega^3(\mathbb{S}^3) \arrow[dlll,out=280,in=80,looseness=0.2,"\Omega^3e",swap] \\
  \Omega^3(\mathbb{S}^1*\mathbb{S}^1) \arrow[r,"\Omega^3\alpha"] &
  \Omega^3(\mathbb{S}^2) \arrow[r,"h"] &
  \Omega^3(\mathbb{S}^1*\mathbb{S}^1) \arrow[r,"\Omega^3(e^{-1})"] &
  \Omega^3(\mathbb{S}^3) \arrow[dlll,out=280,in=80,looseness=0.2,"e_3",swap] \\
  \Omega^2\|\mathbb{S}^2\|_2 \arrow[r,"\Omega\kappa_{2,\mathbb{S}^2}"] &
  \Omega\|\Omega(\mathbb{S}^2)\|_1 \arrow[r,"\kappa_{1,\Omega\mathbb{S}^2}"] &
  \|\Omega^2(\mathbb{S}^2)\|_0 \arrow[r,"e_2"] &
  \Omega(\mathbb{S}^1) \arrow[r,"e_1"] &
  \ZZ
\end{tikzcd}
\]

By implementing this number in a proof assistant with computational
support for univalence and HITs, one should be able to normalise it
using a computer to establish that $\beta = \pm 2$ and hence that
$\pi_4(\bS^3) \cong \ZZ/2\ZZ$. In 2016, by the time Brunerie was
finishing his thesis, there were some experimental proof
assistants based on the cubical type theory of \cite{CCHM18}, but
these were too slow to perform such a complex computation. So, instead
of relying on normalisation, Brunerie spends the second part of
the thesis (Chapters 4--6) to prove, using a lot of the advanced
machinery mentioned above, that $| \beta |$ is propositionally equal to
$2$. However, if one were instead able to compute the
number automatically in a proof assistant, this equality would hold
definitionally---effectively reducing the complexity and length of the
proof by an order of magnitude.

The intriguing possibility of a computer assisted formal proof made
many people interested and countless attempts to normalise Brunerie's
$\beta$ have been made using increasingly powerful computers.
However, to date, no one has succeeded and it is still unclear whether
it is normalisable in a reasonable amount of time. In light of this,
it is natural to wonder whether it is possible to simplify Brunerie's
number in order to be able to compute it. For example, Brunerie's
original definition only involves 1-HITs, as the status of higher HITs
was still quite understudied at the time. With a better understanding
of higher HITs developed in
\cite{LumsdaineShulman19,CoquandHuberMortberg18,CavalloHarper19}, one
quickly sees that the first $3$ maps can be combined into one sending
$1$ to the 3-cell of $\bS^3$ defined as a 3-HIT and not as an iterated
suspension as in Brunerie's thesis. Unfortunately, simple
optimisations like this do not seem to reduce the complexity of the
computation enough and all attempts to run it have thus far failed.

After several unsuccessful attempts at optimising the computation, we
instead decided to formalise the second half of Brunerie's
thesis. However, this is by no means straightforward. The first issue
appears already in the beginning of Chapter~4, a chapter concerning
smash products of spheres. The main result of the section is
Proposition 4.1.2, which says that the smash product is a 1-coherent
symmetric monoidal product on pointed types. However, the proof of
this result is just a sketch and Brunerie writes the following about
it:

\begin{quote}
  \emph{The following result is the main result of this section even
    though we essentially admit it. \cite[Page 90]{Brunerie16}}
\end{quote}

Unfortunately, this result is then used to construct integral
cohomology rings, $H^*(X)$, whose cup product, $\smile$, appears in
the definition of the so-called Hopf invariant which is crucially used
to prove that $| \beta |$ is $2$. While one might be convinced that
Brunerie's informal proof sketch is correct, it is not obvious how one
convinces a proof assistant of this. A complete formalisation would
either have to fill in the holes in the sketch or find an alternative
construction which avoids Proposition 4.1.2. In fact, Brunerie tried
very hard to fill these holes using \Agda
metaprogramming~\cite{brunerie18}. However, he never managed to
typecheck his computer generated proof of the pentagon
identity. Hence, this approach also seems infeasible with current
proof assistant technology.

Luckily, Brunerie, Ljungström and Mörtberg~\cite{BLM22} recently gave
an alternative synthetic definition of the cup product on $H^*(X)$
which completely avoids smash products. This has allowed us to completely skip
the problematic Chapter 4 and, in particular, Proposition 4.1.2, while
still following the proofs in Chapters 5 and 6.  Having a
strategy for a formal proof, we were then able to embark on able to
embark on the ambitious project of formalising Brunerie's proof. Even
though we do not need any theory about smash products, there was still
a lot left to formalise and our final formalisation closely follows
Brunerie's proof, except for various smaller simplifications and
adjustments which we discuss in the paper.

In addition to this, we have also formalised a new proof by
Ljungström~\cite{@Axel-HoTTEST} which completely circumvents Chapters
4--6. This major simplification builds on manually calculating the
image of the element $\eta : \pi_3(\mathbb{S}^2)$, corresponding to
$\beta$ under the isomorphism $\pi_3(\mathbb{S}^2) \cong \ZZ$, by
dividing this isomorphism into several maps, tracing $\eta$ in each
step. In particular, the new proof is completely elementary and does
not rely on advanced tools such as cohomology. The elements that one
obtains while tracing $\eta$ are all new ``Brunerie numbers'' that
should normalise to $\pm 2$. In fact, one of these normalises, in just
under $4$ seconds on a regular laptop, to $-2$ in \CubicalAgda at the
time of writing. So, despite still not being able to compute the
original $\beta$, this work can be seen as an alternative solution to
Brunerie's conjecture about obtaining a computational proof that
$\pi_4(\bS^3) \cong \ZZ/2\ZZ$ which relies on simplifying the Brunerie
number until it becomes effectively computable.

\paragraph{Outline}
The paper closely follows the structure of Brunerie's proof. In
\autoref{sec:hott}, we discuss key results from HoTT that we will need
and their formalisation in \CubicalAgda.
\autoref{sec:homotopygroupsofspheres}, which roughly corresponds to
Chapter 2 of Brunerie's thesis, contains some first results on
homotopy groups of spheres---e.g.\ the computation of
$\pi_n(\mathbb{S}^m)$ for $n \leq m$. We then give Brunerie's
definition of $\beta$ and prove that $\pi_4(\mathbb{S}^3) \cong
\ZZ/\beta\ZZ$, the formalisation of which involves the James
construction and Whitehead products. The remainder of the paper is
then devoted to the formalisation of the different proofs that $\beta
= \pm 2$. We first discuss the formalisation of Chapters 4--6 of
Brunerie's proof in \autoref{sec:brunerieproof}. This involves a lot
of technical machinery like cohomology, the Hopf invariant, etc. We
then, in \autoref{sec:computation}, turn our attention to the new
elementary proof that $\beta = \pm 2$ and the new Brunerie number
which quickly normalises to $-2$ in \CubicalAgda. Here, we also
present some result concerning joins of spheres and the vanishing of
Whitehead products. We conclude in \autoref{sec:conclusion} with a
discussion and comparison of the different formal proofs, as well as
some directions for future work.

Compared to the previous publication on which the current paper is
based, \cite{LICS23}, the main differences are the following.
\begin{itemize}
\item Many proofs which were omitted because of page constraints in
  \cite{LICS23} have been added or extended throughout the paper. In
  particular, the proofs in \cite[Section~VI]{LICS23} have been
  substantially expanded with many details added in
  \autoref{theproof}.
\item In \autoref{sec:computation}, many results from
  \cite[Section~VI]{LICS23} have also been generalised, e.g.\ the
  alternative definition of homotopy groups in terms of
  joins of spheres, $\pi^*_n$, is now studied in general and not just for
  $n = 3$.
\item As part of the expansion and generalisation of
  \cite[Section~VI]{LICS23} in \autoref{sec:computation}, a new
  \autoref{interlude} on joins and smash products of spheres, a new
  \autoref{htpyintermsofjoins} on homotopy groups in terms of joins
  and a new \autoref{standalone} on the possibility of a stand-alone
  proof of Brunerie’s theorem have been added.
\end{itemize}

\paragraph{Formalisation} All results in the paper have been
formalised in \CubicalAgda and are part of the \agdaCubical library,
available at \url{https://github.com/agda/cubical/}. The code in the
paper is mainly literal \Agda code taken verbatim from the library,
but we have taken some liberties when typesetting, e.g.\ shortening
notations and omitting some universe levels. A \CubicalAgda summary
file linking the formalisation and paper can be found at:
\url{https://github.com/agda/cubical/blob/master/Cubical/Papers/Pi4S3-JournalVersion.agda}
The development typechecks with \Agda's \texttt{--safe} flag, which
ensures that there are no admitted goals or postulates.

%% file: sections/hott.tex
\section{Homotopy Type Theory in \CubicalAgda}
\label{sec:hott}

In this section, we concisely summarise the key HoTT concepts
needed for the proofs and their formalisation in \CubicalAgda. This
roughly corresponds to \cite[Chapter~1]{Brunerie16}. For a more
in-depth introduction, see the HoTT Book~\cite{HoTT13} which also
serves as a reference for the formal language ``Book HoTT''. In this
paper, we will present many things with cubical notations, but almost
all of the results also hold with minor changes in Book HoTT where
paths are represented using Martin-Löf's inductive
\textsf{Id}-types~\cite{MartinLof75itt} instead of cubical path
types. In \autoref{sec:conclusion} we discuss in more detail
which proofs crucially rely on cubical features.

All of the results presented in this section were already part of the
\agdaCubical library before we began our formalisation and, while
useful as a resource for our notations, experts on HoTT and
\CubicalAgda can safely skim this section.

\subsection{Elementary HoTT notions and \CubicalAgda notations}

We write $(x : A) \to B\,x$ for dependent function types and denote the
identity function by $\idfun_A : A \to A$. We write
$\Sigma_{x : A} (B\,x)$ for the dependent pair type and \field{fst} and
\field{snd} for its projection maps. In what follows, we mean by a
\emph{pointed type} a dependent pair $(A,\star_A)$ consisting of a
type $A$ and a fixed basepoint $\star_A : A$. For ease of notation, we
will often omit the basepoint and simply write $A$ for the pointed
type $(A,\star_A)$. Given two pointed types $A$ and $B$, the type of
\emph{pointed functions} $A \to_\star B$ consists of pairs
$(f,\star_f)$ where $f : A \to B$ and
$\star_f : \tyPath{f\,\star_A}{\star_B}$ witnesses basepoint
preservation. Again, we simply write $f : A \to_\star B$ and take
$\star_f$ implicit.

HoTT supports inductive types, i.e.\ types inductively generated by
their constructors/points. We write $\func{Bool}$ for the type of
booleans and $\One$ for the unit/singleton type with a single point
$\con{\star_1}$.
A defining feature of HoTT, as opposed to plain Martin-Löf type theory
\cite{MartinLof84bibliopolis}, is the existence of \emph{higher
  inductive types} (HITs). This is a generalisation of inductive types
where we are not only allowed to specify the generating points of the
type in question, but also identifications between these points (and
possibly identifications of these identifications, and so on). This is
useful for defining quotient types, but also for defining spaces when
working in the \emph{types-as-spaces} interpretation of HoTT (see
e.g.\ \cite[Table~1]{HoTT13} and \cite{AwodeyWarren09}). \CubicalAgda
natively supports HITs and a type representing the circle can be
defined as follows:
\ExecuteMetaData[agda/latex/background.tex]{s1}
Here, $\tyPath{\base}{\base}$ denotes the type of identifications of
$\base$ with itself. This is interpreted as the type of \emph{paths}
from $\base$ to itself when viewing $\sphere{1}$ as a space. Hence,
the above HIT captures precisely the representation of the circle as a cell
complex with one $0$-cell ($\base$) and one $1$-cell ($\Loop$). We
always take $\sphere{1}$ to be pointed by $\base$.
In order to discuss the induction principle for $\sphere{1}$, we need
to discuss paths in more detail.
Cubically, paths correspond to functions out of the unit interval,
just like in traditional topology. In \CubicalAgda, there is a
primitive interval type\footnote{For technical reasons, this is
  actually just a ``pre-type'' in \CubicalAgda.} $\func{I}$ with
endpoints \con{i0} and \con{i1}. A path of type $\tyPath{x}{y}$
between two points $x,y:A$ is a function $p : \func{I} \to A$ such
that $p\,\con{i0} = x$ and $p\,\con{i1} = y$ judgmentally. For
instance, $\refl$, the constant path at a point $x$, is defined by:
\ExecuteMetaData[agda/latex/Section2.tex]{refl}
Note that we use ``$=$'' for definitional/judgmental equality and
``$\tyPath{}{}$'' for \CubicalAgda's path-equality. This can be
contrasted with the HoTT Book \cite{HoTT13} which uses the opposite
convention where ``$=$'' is propositional/typal equality and
``$\equiv$'' definitional/judgmental equality.

This type of notational conventions is not the only difference between
\CubicalAgda and Book HoTT. Many proofs that are complicated in Book
HoTT become remarkably direct using the direct treatment of equality
using path types. For instance, function extensionality and its
inverse \func{funExt⁻} are one-liners that just flip the arguments:
\ExecuteMetaData[agda/latex/Section2.tex]{funExt}
\ExecuteMetaData[agda/latex/Section2.tex]{funextInv}
In Book HoTT, however, $\func{funExt}$ is typically proved as a
consequence of the univalence axiom using a rather ingenious proof~\cite{LicataBlog14}
while its inverse follows from path induction. Another elementary
example of a proof involving $\tyPath{\func{\_}}{\func{\_}}$ is
\func{cong} (called \textsf{ap} in Book HoTT), which applies a
function to a path:
\ExecuteMetaData[agda/latex/Section2.tex]{cong}

Although the treatment of paths in \CubicalAgda differs somewhat from
Book HoTT, we may still prove \emph{path induction}: for any dependent
type $B : (y : A)\,(p : \tyPath{x}{y}) \to \func{Type}$, all dependent
functions $f : (y : A)\,(p : \tyPath{x}{y}) \to B\,x\,p$ are uniquely
determined by $f\,x\,(\refl\,x)$. In Book HoTT, this can be used,
among other things, to define the notion of a \emph{dependent} path,
which formalises the situation when two points $a : A$ and $b : B$ are
equal up to a path $p : \tyPath{A}{B}$. In \CubicalAgda, however, the
type of dependent paths is primitive:
\ExecuteMetaData[agda/latex/Section2.tex]{PathP}
In fact, $\tyPath{\func{\_}}{\_}$ is just the special case of
\func{PathP} where the line of paths is constant:
\ExecuteMetaData[agda/latex/Section2.tex]{Path}

We are now ready to describe the induction principle of
$\sphere{1}$. A dependent function $f : (x : \sphere{1}) → B\,x$ is
determined by a point $b : B\,\base$ and a loop
$\ell : \func{PathP}(\lambda i \to B (\Loop\,i))\,b\,b$.  In
\CubicalAgda, this would be written using pattern matching, as in the
left-most definition below, which is introduced side-by-side with the
way it would commonly be written in informal HoTT (as in
Brunerie's thesis):

\begin{minipage}{20cm}
\begin{minipage}[t]{6cm}
\ExecuteMetaData[agda/latex/Background.tex]{s1-fun}
\end{minipage}
\raisebox{.2cm}{
\begin{minipage}[t]{10cm}
\begin{align*}
  &f(\base) = b & \hspace{4.7cm} \\
  &\mathsf{ap}_f(\Loop) = \ell & \hspace{4.7cm}
\end{align*}
\end{minipage}
}
\end{minipage}

\subsection{More higher inductive types}

Let us now introduce the remaining HITs used in
\cite{Brunerie16}. These come equipped with induction principles
analogous to that of $\sphere{1}$. To define higher spheres, we need
suspensions:
\ExecuteMetaData[agda/latex/background.tex]{susp}
We always take suspensions to be pointed by $\north$. We may now
define the $n$-sphere, for $n \geq 1$, by $\sphere{n} =
\SuspN{n-1}{\func{$\mathbb{S}^1$}}$ where $\SuspN{n-1}{}$ denotes
$(n-1)$-fold suspension. We also define $\sphere{-1} = \func{$\bot$}$
(the empty type) and $\sphere{0} = \func{\text{Bool}}$. We remark that
we could equivalently have defined $\sphere{1}$ as the suspension of
$\sphere{0}$ as is done in \cite{Brunerie16}. Our reason for not doing
so is that certain functions using $\sphere{1}$ appear to compute
better with the $\base/\Loop$ definition. Furthermore, this is the
definition used in already existing code in the \agdaCubical library.

We may also capture the (homotopy) pushout of a span $B \xleftarrow{f}
A \xrightarrow{g} C$ by the HIT:
\ExecuteMetaData[agda/latex/background.tex]{pushout}
Diagrammatically this corresponds to:
\[\begin{tikzcd}[ampersand replacement=\&]
	A \& C \\
	B \& {\func{Pushout}\,f\,g}
	\arrow["f"', from=1-1, to=2-1]
	\arrow["g", from=1-1, to=1-2]
	\arrow["\con{inl}",from=2-1, to=2-2,swap]
	\arrow["\con{inr}",from=1-2, to=2-2]
	\arrow["\lrcorner"{anchor=center, pos=0.125, rotate=180}, draw=none, from=2-2, to=1-1]
\end{tikzcd}\]

We use pushouts to define the wedge sum of two pointed types, denoted
${A \,\func{$\vee$}\, B}$, the join of two types, denoted
$\JJoin{A\,}{B}$, and the cofibre of a map $f :A \to B$, denoted
$\cofib{f}$:
\[\begin{tikzcd}[ampersand replacement=\&]
	\One\& B \& {A \,\func{$\times$}\, B} \& B \& A \& B\\
	A \& {A \,\func{$\vee$}\, B} \& A \& {{A}\,\func{\text{\textasteriskcentered}}{B}} \& \One \& \cofib{f}
	\arrow[from=1-1, to=2-1]
	\arrow[from=1-1, to=1-2]
	\arrow[from=2-1, to=2-2]
	\arrow[from=1-2, to=2-2]
	\arrow["\field{fst}"', from=1-3, to=2-3]
	\arrow["\field{snd}", from=1-3, to=1-4]
	\arrow[from=2-3, to=2-4]
	\arrow[from=1-4, to=2-4]
        \arrow["f",from=1-5, to=1-6]
        \arrow[from=1-5, to=2-5]
        \arrow[from=1-6, to=2-6]
        \arrow[from=2-5, to=2-6]
        \arrow["\lrcorner"{anchor=center, pos=0.125, rotate=180}, draw=none, from=2-6, to=1-5]
        \arrow["\lrcorner"{anchor=center, pos=0.125, rotate=180}, draw=none, from=2-4, to=1-3]
        \arrow["\lrcorner"{anchor=center, pos=0.125, rotate=180}, draw=none, from=2-2, to=1-1]
\end{tikzcd}\]

Two particularly important functions out of wedge sums are
\ExecuteMetaData[agda/latex/background.tex]{nabla}
and
\ExecuteMetaData[agda/latex/background.tex]{totimes}

\subsection{Truncation levels and $n$-truncations}

An important concept in HoTT is that of Voevodsky's
\emph{h-levels}~\cite{Voevodsky10bonn}, which gives rise to the notion
of an \emph{n-type}. Since types in HoTT are interpreted as spaces (or
rather, as homotopy types), they are not only determined by their
points but also by which higher paths they may contain.  We say that a
type $A$ is an $n$-type if all $(n+1)$-dimensional structure of $A$ is
trivial.  Formally, this is captured by an inductive definition. We
say that $A$ is a $(-2)$-type if it is contractible, i.e.\ consisting
of a single point, as captured by $\func{isContr}\,A =
\Sigma_{a_0 : A}((a : A) \to \tyPath{a_0}{a})$.
We inductively say that $A$ is an $(n+1)$-type if for any $x ,y
: A$, the type $\tyPath{x}{y}$ is an $n$-type. We call $(-1)$-types
\emph{propositions} and $0$-types \emph{sets}.

We can turn any type $A$ into an $n$-type by $n$-\emph{truncation},
denoted $\truncT{n}{A}$. For instance, the $(-1)$-truncation may be
directly defined using the following HIT:
\ExecuteMetaData[agda/latex/background.tex]{ptrunc}
We often use direct definitions like this of $(-1)$- and
$0$-truncation in our formalisation, and similar constructions work
for any fixed value of $n$, but not when $n$ is arbitrary. For higher
$n$ we rely on the hub-and-spoke
construction~\cite[Section~7.3]{HoTT13}.
\ExecuteMetaData[agda/latex/background.tex]{trunc}

One caveat with truncations is that a map $f : A \to B$ does \emph{not}, in
general, induce a map $f : \truncT{n}{A} \to B$. This is, however, the
case when $B$ is an $n$-type. In particular, $f$ always induces a
function $\truncT{n}{f} : \truncT{n}{A} \to \truncT{n}{B}$.

\subsection{Univalence, loop spaces, and H-spaces}

In order to introduce Voevodsky's univalence
principle~\cite{Voevodsky10cmu}, we need to define the (homotopy)
fibre of a function. Given a function $f : A \to B$ and a point $b :
B$, we define the fibre of $f$ over $b$ by $\fib{f}\,b = \Sigma_{x :
  A}{(\tyPath{f\,a}{b}})$. We say that $f : A \to B$ is an
equivalence, written $f : \tyEquiv{A}{B}$, if $\fib{f}\,b$ is
contractible for all $b : B$. In order to prove that a function $f : A
\to B$ is an equivalence, it suffices to provide an inverse $f^- : B
\to A$ and two paths $\tyPath{f\circ f^-}{\idfun_B}$ and
$\tyPath{f^-\circ f}{\idfun_A}$.  If $f$ is also pointed, we write $f
: A \,\func{$\simeq_\star$}\,B$.

Univalence states that the canonical map $ \tyPath{A}{B} \to
\tyEquiv{A}{B}$, defined by path induction, is an equivalence. In
particular, we get a map $\func{ua} : \tyEquiv{A}{B} \to
\tyPath{A}{B}$ promoting equivalences to paths. This provides us with
a useful method for transferring proofs between equivalent types which
extends to \emph{structured} types and are then referred to as the
\emph{structure identity principle}~\cite[Section 9.8]{HoTT13}.

Transferring proofs is, however, not the only use case of univalence
in HoTT. It can also be used to characterise \emph{loop spaces} of
HITs. This is often done using the \emph{encode-decode
  method}~\cite[Section~8.1.4]{HoTT13}, a type theoretic analogue of
proofs by contractibilty of total spaces of fibrations. In HoTT, we
define the loop space of a pointed type $A$, by $\LoopSpace{A} =
(\tyPath{\star_A}{\star_A})$. This is again pointed by
$\refl\,\star_A$, so we may iterate this definition to get the $n$th
loop space of $A$, denoted $\LoopSpaceN{n}{A}$. Loop spaces belong to
a particularly important class of types called \emph{H-spaces}. These consist
of a pointed type $B$ equipped with a unital magma structure
\begin{align*}
&\mu : B \,\func{$\times$} B \to B \\
&\mu_l : (b : B) \to \tyPath{\mu(\star_B,b)}{b} \\
&\mu_r : (b : B) \to \tyPath{\mu(b,\star_B)}{b}
\end{align*}
satisfying $\tyPath{\mu_l\,\star_B}{\mu_r\,\star_B}$.  Another
particularly important H-space for our purposes is $\sphere{1}$, for
which we will use $\func{$+$}$ to denote its binary
operation. $\sphere{1}$ also comes equipped with a notion of inversion
which we will denote by $\func{$-$}$. In fact, $\sphere{1}$ is a
commutative and associative H-space.

%% file: sections/homotopygroupsofspheres.tex
\section{First results on homotopy groups of spheres}
\label{sec:homotopygroupsofspheres}

In this section, we cover \cite[Chapter~2]{Brunerie16}, which
introduces some elementary results on the homotopy groups of
spheres. All of these results can also be found in the HoTT
Book~\cite{HoTT13}.  Before even stating them, we need homotopy
groups:
\begin{defi}[Homotopy groups]\label{def:hom-gr}
  For $n : \NN$, we define the $n$th \emph{homotopy group} of a
  pointed type $A$ by:
  \[\func{$\pi$}_n(A) = \truncT{0}{\sphere{n} \to_\star A} \]
\end{defi}
The name homotopy \emph{group} should be taken with a grain of salt:
it, in general, only has a group structure when $n \geq 1$ (abelian
when $n \geq 2$). The structure may be defined, much like in \cite[Section
  5]{BuchholtzFavonia18}, by considering the
equivalence
$\tyEquiv{(\sphere{n} \to_\star A)}{(\sphere{n-1}
  \to_\star\LoopSpace{A})}$, where the latter type has a
multiplication given by pointwise path composition.
An alternative definition of $\func{$\pi$}_n(A)$ is via loop
spaces. There is an equivalence
$\omega_n : \tyEquiv{\LoopSpaceN{n}{A}}{(\sphere{n} \to_\star A)}$
and, hence, we could equivalently have defined $\func{$\pi$}_n(A)$ by
setting $\func{$\pi$}_n(A) = \truncT{0}{\LoopSpaceN{n}{A}}$.  This
makes the group structure on $\func{$\pi$}_n(A)$ more transparent: it
is simply path composition. This is the definition used in the HoTT
Book~\cite{HoTT13}. Brunerie uses both definitions in his thesis and
often passes between the two without comment.

An elementary but crucial result for the computation of homotopy
groups is the existence of the \emph{long exact sequence of homotopy
  groups}. Its proof is usually phrased using the loop space
definition of homotopy groups as in e.g.\ \cite[Theorem
8.4.6]{HoTT13}.  For ease of notation, let us write $\fib{f}$ for the
fibre of a pointed function $f : A \to_\star B$ over the basepoint of
$B$.
\begin{prop}[LES of homotopy groups]\label{prop:LES}
  For any pointed map $f : A \to_\star B$, there is a long exact
  sequence
  \textnormal{
\[\begin{tikzcd}[ampersand replacement=\&,row sep = -2pt, column sep = 38pt]
	\& \dots \& {\func{$\pi$}_{n+1}(B)} \\
	\& \phantom{.} \& \\
	{\func{$\pi$}_{n}(\fib{f})} \& {\func{$\pi$}_n(A)} \& {\func{$\pi$}_n(B)} \\
	\& \phantom{.} \& \\
	{\func{$\pi$}_{n-1}(\fib{f})}  \& \dots
	\arrow[from=3-1, to=3-2]
	\arrow[from=3-2, to=3-3]
	\arrow[from=1-2, to=1-3]
	\arrow[from=2-2, to=3-1,shorten <=-7.5pt, curve={height=6pt}]
        \arrow[from=4-2, to=5-1,shorten <=-7.5pt, curve={height=6pt}]
	\arrow[shorten <=4pt,shorten >=-5pt, no head, from=1-3, to=2-2,curve={height=-3pt}]
	\arrow[shorten <=4pt,shorten >=-5pt, no head, from=3-3, to=4-2,curve={height=-3.5pt}]
	\arrow[from=5-1, to=5-2]
\end{tikzcd}\]}
where the horizontal maps are induced by the functorial action of $\func{$\pi$}_n$ on \textnormal{$\field{\mathsf{fst}} : \fib{f} \to A$} and $f : A \to B$. 
\end{prop}
Above, we have implicitly taken the kernel and image of a group homomorphism $\phi : G \to H$ to be defined by
\begin{align*}
  \func{ker}\,\phi\, &= \func{fib}\,\phi\,0_H\\
  \func{im}\,\phi &= \Sigma_{h : H} \truncT{-1}{\Sigma_{g : G} (\phi(g) \,\func{$\equiv$}\,h)}
\end{align*}

When analysing loop spaces and homotopy groups of suspensions, the
following function is of great importance. It will be used in many
constructions to come.
\begin{defi}[The suspension map]
  Given a pointed type $A$, there is a canonical map
  $\func{$\sigma$}: A \to \LoopSpace (\Susp{A})$ given by
  \[
    \func{$\sigma$}\,x = \merid{x} \,\func{$\cdot$}\, \sym{(\merid{\star_A})}
  \]
\end{defi}
This induces a homomorphism on homotopy groups by post-composition:
\[
  \func{$\pi$}_n(A)
  \xrightarrow{\,\func{$\sigma$}_*\,}
  \func{$\pi$}_n(\LoopSpace{(\Susp{A})})
  \xrightarrow{\,\func{$\cong$}\,}
  \func{$\pi$}_{n+1}(\Susp{A})
\]
We will often, with some abuse of notation, simply write $\func{$\sigma$}_*$ for this composition.
We also define
$\func{$\sigma$}_n : \truncT{n}{A} \to
\LoopSpace{\truncT{n+1}{\Susp{A}}}$ by \textnormal{
    \[\func{$\sigma$}_n\trunc{x} = \ap{\trunc{\con{\!\_\!}}}{(\func{$\sigma$}\,x)}\]
}%
We will soon see the suspension map in action, but first we need the
following elementary result.
\begin{prop}[Join of spheres]\label{prop:spherejoin}
  $\tyEquiv{\JJoin{\sphere{n}}{\,\sphere{m}}}{\sphere{n+m+1}}$.
\end{prop}
In fact, as we will see in~\autoref{sec:computation}, there is more to
say about this equivalence. We make a forwards reference to \autoref{prop:F-equiv} and the preceding discussion for a detailed account of its construction.

In particular, \autoref{prop:spherejoin} gives us an equivalence
$\tyEquiv{\JoinS}{\sphere{3}}$.
Using this fact, we define the following map, which will play a
crucial role in the analysis of $\func{$\pi$}_4(\sphere{3})$.
\begin{defi}[Hopf map]\label{def:hopfmap}
  We define $\hopf : \sphere{3} \to \sphere{2}$ by the composition
  $\sphere{3} \xrightarrow{\sim} \JoinS \xrightarrow{\prehopf}
  \sphere{2}$ where $\prehopf$ is given by \textnormal{
  \ExecuteMetaData[agda/latex/BruneriesProof.tex]{hopf}
  }%
where $\func{$-$}$ is defined using the H-space and inversion structure on $\sphere{1}$.
\end{defi}
It turns out that the following is true~\cite[Theorem~8.5.1]{HoTT13}.
\begin{prop}[The fibre of the Hopf map]\label{prop:fib-hopf}
  The fibre of $\hopf$ is equivalent to $\sphere{1}$, i.e.\
  $\tyEquiv{\fib{\hopf}}{\sphere{1}}$.
\end{prop}
\autoref{prop:fib-hopf} gives us a fibration sequence $\sphere{1} \to
\sphere{3} \to \sphere{2}$ which, in particular, will allow us to
connect homotopy groups of $\sphere{2}$ with those of $\sphere{3}$ and
$\sphere{1}$.
For this, we need to introduce the notion of \emph{connectedness}. We
say that a type $A$ is $n$-connected if $\truncT{n}{A}$ is
contractible. Similarly, we say that a function $f : A \to B$ is
$n$-connected if all of its fibres are $n$-connected. This means, in
particular, that the induced function
$\truncT{n}{f} : \truncT{n}{A} \to \truncT{n}{B}$ is an equivalence.
The following is an immediate consequence of the definition of
$n$-truncations.
\begin{lem}[Connectedness of spheres]\label{lem:sphere-conn}
  For $n \geq -1$, $\sphere{n}$ is $(n-1)$-connected.
\end{lem}
Using~\autoref{lem:sphere-conn}, we can easily prove
the following:
\begin{prop}[{\cite[Proposition~2.4.1]{Brunerie16}}]\label{prop:low-sphere}
  For $n < m$, the group $\func{$\pi$}_n(\sphere{m})$ is trivial.
\end{prop}
For the sake of completeness, let us take the liberty of mentioning some
results from \cite[Chapter~3]{Brunerie16} already here, since they
also concern low-dimensional homotopy groups of spheres. A crucial
result is the following theorem~\cite[Theorem~8.6.4]{HoTT13}:
\begin{thm}[Freudenthal suspension theorem]\label{thm:Freudenthal}
  Given an $n$-connected and pointed type $A$, the map
  $\func{$\sigma$}: A \to \LoopSpace (\Susp{A})$ is $2n$-connected.
\end{thm}
On can easily deduce from~\autoref{thm:Freudenthal} that, in particular,
  $\func{$\sigma$}_n : \truncT{n}{A} \to \truncT{n}{\LoopSpace (\Susp{A})}$
  is an equivalence. This allows us to prove the following result:
\begin{cor}\label{cor:stable}
  For $n \geq 1$, we have $\tyIso{\func{$\pi$}_n(\sphere{n})}{\bZ}$. Furthermore,
  $\func{$\pi$}_n(\sphere{n})$ is generated by $i_n =
  \trunc{\idfun_{\sphere{n}}}$.
\end{cor}
\begin{proof}
  The synthetic proof of the classical result that
  $\tyIso{\func{$\pi$}_1(\sphere{1})}{\bZ}$ is due to Licata and
  Shulman~\cite{LicataShulman13}.  The fact that
  $\tyIso{\func{$\pi$}_2(\sphere{2})}{\func{$\pi$}_1(\sphere{1})}$ is
  given by the LES associated to the Hopf fibration combined
  with~\autoref{prop:low-sphere}. The fact that
  $\tyIso{\func{$\pi$}_{n+1}(\sphere{n+1})}{\func{$\pi$}_n(\sphere{n})}$
  is an immediate consequence of~\autoref{thm:Freudenthal}. The second
  statement follows by induction on $n$, using that suspension is
  functorial and thereby preserves the identity map.
\end{proof}

We have now analysed all homotopy groups $\func{$\pi$}_n(\sphere{m})$ with
$n \leq m$. This yields the following:
\begin{prop}\label{prop:hopfcomp}
  Post-composition by \textnormal{$\hopf$} induces an isomorphism
  $\tyIso{\func{$\pi$}_3(\sphere{3})}{\func{$\pi$}_3(\sphere{2})}$.
\end{prop}
\begin{proof}
  By \autoref{prop:LES} and \autoref{prop:fib-hopf}, we get an exact sequence
  \[\func{$\pi$}_3(\sphere{1}) \to \func{$\pi$}_3(\sphere{3}) \xrightarrow{\hopf_*} \func{$\pi$}_3(\sphere{2}) \to \func{$\pi$}_2(\sphere{1})\]
  as $\func{$\pi$}_n(\sphere{1})$ vanishes for $n> 1$, $\hopf_*$
  is an isomorphism.
\end{proof}
\begin{cor}\label{prop:pi3-gen}
  There is an isomorphism
  $\func{$\psi$}:\tyIso{\func{$\pi$}_3(\sphere{2})}{\bZ}$. Furthermore,
  $\func{$\pi$}_3(\sphere{2})$ is generated by $\hopf$.
\end{cor}
\begin{proof}
  By \autoref{cor:stable} we know that $\func{$\pi$}_3(\sphere{3})$ is
  generated by the identity function on $\sphere{3}$. We know that the
  isomorphism
  $\tyIso{\func{$\pi$}_3(\sphere{3})}{\func{$\pi$}_3(\sphere{2})}$ is
  given by post-composition by $\func{hopf}$ and thus the generator of
  $\func{$\pi$}_3(\sphere{3})$ is mapped to $\func{hopf}$.
\end{proof}

\subsection{Formalisation of Brunerie's Chapter 2}
Most of these results have already been added to \agdaCubical by
Mörtberg~\&~Pujet~\cite{cubicalsynthetic},
Ljungström~\cite{LjungstromMsc20}, and Brunerie,
Ljungström \& Mörtberg~\cite{BLM22}.
The Freudenthal suspension theorem was formalised in \CubicalAgda by
Cavallo~\cite{FreudenthalFormalization}, using a direct cubical
proof following~\cite[Theorem~8.6.4]{HoTT13}.  \autoref{cor:stable}
was given a direct proof, following the computation of cohomology
groups of spheres in~\cite{BLM22}.

There were some technical difficulties related to the equivalence
$\omega_n : \tyEquiv{\LoopSpaceN{n}{A}}{ (\sphere{n}\to_\star A)}$,
which is used to show that the two different definitions of homotopy
groups are equivalent. In several proofs, it is more natural to work
on the left-hand-side of $\omega_n$. At the same time, working on the
right-hand-side often makes constructing elements easier (compare, for
instance, an explicit description of the generator of
$i_3 : \func{$\pi$}_3(\sphere{3})$ described as a $3$-loop in
$\sphere{3}$ to the very compact definition
$i_3 = \trunc{\idfun_{\sphere{3}}}$). This means that we often have to
translate between the two definitions. One particularly important
example is the LES of homotopy groups associated to a function
$A \to_\star B$. On each level, the maps are given as follows:
  \begin{align*}
    \LoopSpaceN{n}{(\fib{f})}
    \xrightarrow{\,\LoopSpaceN{n}{\mathsf{fst}}\,} \LoopSpaceN{n}{A}
    \xrightarrow{\,\LoopSpaceN{n}{f}\,} \LoopSpaceN{n}{B}
  \end{align*}
  This is then transported to the definition of homotopy groups as
  maps from spheres via $\omega_n$. For the proof of e.g.\
  \autoref{prop:pi3-gen}, we need to know that the maps in the
  sequence are given as follows:
  \begin{align*}
    \func{$\pi$}_n{(\fib{f})}
    \xrightarrow{\,\field{fst}_*\,} \func{$\pi$}_{n}{(A)}
    \xrightarrow{\,f_*\,} \func{$\pi$}_{n}{(B)}
  \end{align*}
  What we need is then more than just an equivalence $\omega_n :
  \tyEquiv{\LoopSpaceN{n}{A}}{(\sphere{n}\to_\star A)}$ -- we need to
  show that this equivalence is functorial. This is implicitly assumed
  in Brunerie's thesis, but, in \CubicalAgda, we need to make it
  precise. Formalising this fact is not entirely trivial. First, we
  need a tractable definition of the equivalence in question. It can
  be described inductively with base case $\omega_1 : \LoopSpace{A}
  \to (\sphere{1} \to_\star A)$ given by:
\begin{align*}
  \omega_1\,p\,\base &= \star_A\\
  \omega_1\,p\,(\Loop\,i) &= p\,i
\end{align*}
which we take to be pointed by $\refl$. It is easy to verify that this
is an equivalence. We define $\omega_{n+1}$ by the composition:
\begin{align*}
  \LoopSpaceN{n+1}{A} = \LoopSpace{(\LoopSpaceN{n}{A})}
  \xrightarrow{\,\LoopSpace{\omega_n}\,}
  \LoopSpace{(\sphere{n} \to_\star A)}
  \xrightarrow{\func{funExt$^-_\star$}}
  (\sphere{n} \to_\star \LoopSpace A)
  \xrightarrow{\,}
  (\sphere{n+1} \to_\star A)
\end{align*}
where the last arrow comes from the adjunction
$\func{Susp} \dashv \LoopSpace{}$. This is a composition of
equivalences, and hence an equivalence. We then need to verify that
the following commutes
\[\begin{tikzcd}[ampersand replacement=\&]
	{\LoopSpaceN{n}{A}} \&\& {(\sphere{n} \to_\star A)} \\
	\\
	{\LoopSpaceN{n}{B}} \&\& {(\sphere{n} \to_\star B)}
	\arrow["{\omega_{n}}", from=1-1, to=1-3]
	\arrow["{\LoopSpaceN{n}{\!f}}"', from=1-1, to=3-1]
	\arrow["{\omega_{n}}"', from=3-1, to=3-3]
	\arrow["{f_*}", from=1-3, to=3-3]
      \end{tikzcd}\]
This can be proved inductively. The base case is easy and the
inductive step is given by the following diagram
\[\begin{tikzcd}[ampersand replacement=\&]
	\&\&\& {\LoopSpace(\sphere{n} \to_\star A)} \\
	{\LoopSpaceN{n+1}{A}} \&\& {(\sphere{n+1} \to_\star A)} \\
	\\
	{\LoopSpaceN{n+1}{B}} \&\& {(\sphere{n+1} \to_\star B)} \\
	\&\&\& {\LoopSpace(\sphere{n} \to_\star B)}
	\arrow["{\omega_{n+1}}"', from=2-1, to=2-3]
	\arrow["{\LoopSpaceN{n+1}{f}}"', from=2-1, to=4-1]
	\arrow["{\omega_{n+1}}", from=4-1, to=4-3]
	\arrow["{f_*}", from=2-3, to=4-3]
	\arrow["\tyEquiv{}{}", from=1-4, to=2-3]
	\arrow["{\LoopSpace{f_*}}", from=1-4, to=5-4]
	\arrow["\tyEquiv{}{}"', from=5-4, to=4-3]
	\arrow["{\LoopSpace{\omega_n}}"', from=4-1, to=5-4]
	\arrow["{\LoopSpace{\omega_n}}", from=2-1, to=1-4]
\end{tikzcd}\]
where the commutativity of the outer square comes from the base case
paired with the inductive hypothesis, the triangles from the
definition of $\omega_{n+1}$ and the right-most square from a
straightforward argument.

%% file: sections/brunerienumber.tex
\section{The Brunerie number}
\label{sec:brunerienumber}

Here we give an overview of the first half of Brunerie's proof. This
corresponds to \cite[Chapter~3]{Brunerie16} and culminates in the
isomorphism $\tyIso{\func{$\pi$}_4(\sphere{3})}{\bZ/\brunerie\bZ}$ for
an at this point unknown ``Brunerie number'' $\brunerie : \bZ$. We
also discuss the formalisation of this part of the proof and various
simplifications found during the formalisation.

\subsection{The James construction}

To define $\brunerie$, Brunerie uses the \emph{James
  construction}~\cite{James55}, which he introduced in HoTT and
partially formalised in \cite{GB-James}.
\begin{prop}[James construction]\label{prop:James}
  For a $(k\geq 0)$-connected pointed type $A$, there are types
  $\james{n}{A}$ with inclusions
  \[
    \james{0}{A} \lhook\joinrel
    \xrightarrow{\,j_0\,} \james{1}{A}
    \lhook\joinrel
    \xrightarrow{\,j_1\,} \james{2}{A}
    \lhook\joinrel
    \xrightarrow{\,j_2\,} \cdots
    \]
    such that its sequential colimit
    $\tyEquiv{\james{\infty}{A}}{\LoopSpace{(\Susp{A})}}$. Furthermore,
    $j_n : \james{n}{A} \hookrightarrow \james{n+1}{A}$ is
    $(n(k+1)+(k-1))$-connected.
  \end{prop}
A consequence of \autoref{prop:James} is the following fact
\begin{prop}\label{cor:J2-conn}
  Given a $(k \geq 0)$-connected type $A$, there is a
  $(3k+1)$-connected map $\james{2}{A} \xrightarrow{} \LoopSpace{(\Susp{A})}$.
\end{prop}
The proof of~\autoref{cor:J2-conn} uses that $\james{\infty}{A}$, the
sequential colimit of the sequence in~\autoref{prop:James}, can be
shown to be equivalent to $\LoopSpace{(\Susp{A})}$. This, paired with
some results on the connectivity of sequential colimits, gives the
statement. A key consequence of this is the following result which allows us to express ${\func{$\pi$}_4(\sphere{3})}$ as ${\func{$\pi$}_3(\js)}$ -- a group which turns out to be quite a bit easier to reason about.

\begin{thm}\label{cor:P4S3-J2S2}
  $\tyIso{\func{$\pi$}_4(\sphere{3})}{\func{$\pi$}_3(\js)}$
\end{thm}
\begin{proof}
Because $\sphere{2}$ is $1$-connected, \autoref{cor:J2-conn} tells us
that there is a $4$-connected map
\[
\js \to \LoopSpace(\Susp{\sphere{2}}) = \LoopSpace (\sphere{3})
\]
  In particular, it is $3$-connected and induces an equivalence
  $\tyEquiv{\truncT{3}{\js}}{\truncT{3}{\LoopSpace{\sphere{3}}}}$. We
  get:
\[
  \func{$\pi$}_4(\sphere{3}) \,\,\func{$\cong$}\,\, \func{$\pi$}_3(\LoopSpace{\sphere{3}}) \,\,\func{$\cong$}\,\,
  \func{$\pi$}_3(\js) \qedhere
\]
\end{proof}

\subsection{Formalisation of the James construction}
This is a particularly technical part of Brunerie's thesis, primarily
due to the high number of higher coherences which need to be verified
in the proof of \autoref{prop:James}. While this has, subsequent to our
efforts, been formalised in its entirety by Kang~\cite{JamesFormalization}, we
have taken a shortcut by giving a direct proof of~\autoref{cor:P4S3-J2S2}, which means
we do not in fact need the full James construction.
Consequently, we instead give direct definitions of $\james{n}{A}$ for $n \leq 2$
for a pointed type $A$.
\begin{defi}[Low dimensional James construction]\label{def:J1-J2}
  We define \textnormal{$\james{0}{A} = \One$} and $\james{1}{A} = A$.
  The type $\james{2}{A}$ is defined as the pushout:
  \[
    \begin{tikzcd}[ampersand replacement=\&]
        {A \vee A} \& A \\
        {A \times A} \& {\james{2}{A}}
        \arrow["{\ivee}"', from=1-1, to=2-1]
        \arrow["{\Nabla}", from=1-1, to=1-2]
        \arrow[from=2-1, to=2-2]
        \arrow[from=1-2, to=2-2]
        \arrow["\lrcorner"{anchor=center, pos=0.125, rotate=180}, draw=none, from=2-2, to=1-1]
   \end{tikzcd}
    \]
\end{defi}
We remark that the construction in~\autoref{def:J1-J2} is not definitionally the same as Brunerie's; in his thesis, these constructions are
theorems rather than definitions. Here we take them as
definitions. With $\james{n}{A}$ defined this way, the map $j_0 : \james{0}{A} \to \james{1}{A}$
is just the constant pointed map and $j_1 : \james{1}{A}
\to \james{2}{A}$ is $\con{inr}$.

Before we continue, let us temporarily redefine $\sphere{2}$ to be the following
equivalent HIT. This will make some of the following constructions more
compact.
\ExecuteMetaData[agda/latex/background.tex]{s2}

The next lemma will be crucial. It is a special case of the \emph{Wedge Connectivity
  Lemma}~\cite[Lemma~8.6.2]{HoTT13}, of which we have formalised a version
of the proof of the sphere case in~\cite[Lemma~8]{BLM22}. From the
point of view of formalisation, this proof is easier to work with
since it gives more useful definitional equalities.
\begin{lem}[Wedge connectivity for $\sphere{2}$]\label{lem:wedge-conn}
Let $P : \sphere{2} \times \sphere{2} \to \nType{2}$. Any
function $f : (x : \sphere{2}\,\func{$\times$}\, \sphere{2}) \to P\,x$
is induced by the following data:
\begin{align*}
    f_l &:(x: \sphere{2}) \to P(x\,\con{,}\,\base)\\
    f_r &:(y: \sphere{2}) \to P(\base\,\con{,}\,y)\\
    f_{lr} &: \tyPath{f_l\,\base}{f_r\,\base}
\end{align*}
\end{lem}
Before we discuss the formalisation of \autoref{cor:P4S3-J2S2} stated
with the low dimensional James construction, we first construct the
following function. The goal is to define a family of equivalences
$f_x : \tyEquiv{\JTrunc}{\JTrunc}$ over $x : \sphere{2}$. We do this
by truncation elimination and pattern matching on $x$, starting with
the $\base$ case:
\begin{align*}
  f_{\base}\,\trunc{\inl(x,y)} &= \trunc{\inl(x\,\con{,}\,y)} \\
  f_{\base}\,\trunc{\inr{z}} &= \trunc{\inl(\base\,\con{,}\,z)} \\
  f_{\base}\,\trunc{\push{(\inl{x}})\,i} &= \trunc{{(\push{(\inl{x})} \,\func{$\cdot$}\,\sym{\push{(\inr{x})}})}\,i} \\
  f_{\base}\,\trunc{\push{(\inr{y}})\,i} &= \trunc{\inl(\base\,\con{,}\,y)} \\
  f_{\base}\,\trunc{\push{(\push{y}\,j)}\,i} &= \dots
\end{align*}
where the omitted step consists of a proof that
$\push{(\inl{\base})} \cdot \sym{\push{(\inr{\base})}}\,
\func{$\equiv$} \,\refl$. It is an easy lemma that $f_\base$ is equal
to the identity on $\JTrunc$. To complete the definition of $f_x$, we
need to consider the case when $x = \surf\,i\,j$. This amounts to
providing a dependent function:
\[
  f_\surf : (x : \JTrunc) \to \LoopSpaceN{2}{(\JTrunc\,\con{,}\,f_\base\,x)}
\]
To do this, we will, in particular, need to provide a family of fillers
\[Q_{(x\,\con{,}\,y)} : \refl_{\trunc{\inl(x\,\con{,}\,y)}}
\,\func{$\equiv$}\,\refl_{\trunc{\inl(x\,\con{,}\,y)}}\]
This is a $1$-type, and
thus~\autoref{lem:wedge-conn} applies. We define:
\begin{align*}
  Q_{(\base\,\con{,}\,y)}\,i\,j &= \trunc{\inl(\surf\,i\,j\,\con{,}\,y)} \\
  Q_{(x\,\con{,}\,\base)} \,i\,j &= \trunc{\inl(x\,\con{,}\,\surf\,i\,j)}
\end{align*}
The fact that these two constructions agree when both $x$ and $y$ are
$\base$ is a technical but relatively straightforward lemma. Thereby,
$Q_{(x\,\con{,}\,y)}$ is defined. We may now define $f_\surf$:
\begin{align*}
  f_\surf \,\trunc{\inl(x\,\con{,}\,y)} &= Q_{(x\,\con{,}\,y)}\\
  f_\surf \,\trunc{\inr{z}} &= Q_{(\base\,\con{,}\,z)}
\end{align*}
The higher cases are easy due to the fact that the goal becomes
$0$-truncated, making it sufficient to define them for
$\base : \sphere{2}$. Thus, $f_x$ is defined for all $x : \sphere{2}$.
\begin{lem}
  For $x : \sphere{2}$, $f_x$ is an automorphism on $\JTrunc$.
\end{lem}
\begin{proof}
  To make coming proofs easier, this is proved by explicitly
  constructing the inverse analogously to $f_x$.
  \begin{align*}
    f^{-1}_\base\,x &= f_\base\,x \\
    f^{-1}_{\surf}\, x &= f_{\sym{\surf}}\,x
  \end{align*}
  Proving that these cancel is technical, but direct.
\end{proof}
We are now ready to prove the following statement, which is a rephrasing
of \autoref{cor:P4S3-J2S2}.
\begin{prop}\label{prop:J2-direct}
  $\LoopSpace{\truncT{4}{\sphere{3}}}\,\func{$\simeq$}\, \truncT{3}{\james{2}{\sphere{2}}}$
\end{prop}
\begin{proof}
  We take $\sphere{3} = \Susp{\sphere{2}}$,
  where $\sphere{2}$ is defined using $\base$/$\surf$ as above. We
  employ the encode-decode method and define a family of $3$-types
  over $\truncT{4}{\sphere{3}}$. Since the universe of $3$-types is a
  $4$-type, we may do so by truncation elimination:
  \begin{align*}
    \Code &: \truncT{4}{\sphere{3}} \to \textsf{$3$\func{-Type}} \\
    \Code\,\trunc{\north}\, &= \JTrunc \\
    \Code\,\trunc{\south} &= \JTrunc \\
    \Code\,\trunc{\merid{x}\,i} &= \func{ua}\,f_x\,i
  \end{align*}
  We now need to define two families of functions
  \begin{align*}
  &\encode_{x} : \trunc{\north}
  \,\func{$\equiv$}\,x \to \Code\,x \\
  &\decode_{x} :
  \Code\,x \to \trunc{\north} \,\func{$\equiv$}\,x
  \end{align*}
  over $x :
  \truncT{4}{\sphere{3}}$. We define $\encode_x$ by path induction, sending $\refl$ to the basepoint in $\JTrunc$. We define $\decode_x$ by
  truncation elimination and pattern matching on $x$. The crucial step
  is defining $\decode_{\trunc{\north}} : \JTrunc
  \to \LoopSpace{\truncT{4}{\sphere{3}}}$. On point constructors, it
  is given by
  \begin{align*}
    \decode_{\trunc{\north}}\,(\inl{(x\,\con{,}\,y)}) &= \func{$\sigma$}\,x
    \,\func{$\cdot$} \, \func{$\sigma$} \,y
    \\ \decode_{\trunc{\north}}\,(\inr{z}) &= \func{$\sigma$}\,z
  \end{align*}
  which is easily verified to be coherent with the higher
  constructors.
  The case $\decode_{\trunc{\south}}$ is immediately induced by
  $\decode_{\trunc{\north}}$, since $\tyPath{\north}{\south}$ via
  $\merid{\base}$.  The case $\decode_{\trunc{\merid{a}\,i}}\,y$
  amounts to showing that
  \[
    \decode_{\trunc{\north}}(f^{-1}_a\,y) \,\,\func{$\equiv$}\,\,\decode_{\trunc{\north}}\,y \,\func{$\cdot$} \,
  \sym{(\func{$\sigma$}\,\trunc{a})}
\]
The proof is technical but is
  greatly aided by~\autoref{lem:wedge-conn}. The fact that
  $\decode_{x}(\encode_{x}\,p) \, \func{$\equiv$}\, p$ for each $p :
  \north\, \func{$\equiv$} \,x$ holds by path induction. Finally, the
  fact that $\encode_{\north}(\decode_\north\,y) \func{$\equiv$} \, y$
  holds for each $y : \JTrunc$ holds by some technical but simple
  path algebra. Hence $\decode_{\trunc{\north}}
  : \LoopSpace{\truncT{4}{\sphere{3}}} \to \JTrunc$ is an equivalence.
\end{proof}
We get \autoref{cor:P4S3-J2S2} as an immediate corollary of
\autoref{prop:J2-direct} via the same sequence of isomorphisms as in the
proof of~\autoref{cor:P4S3-J2S2}.

\subsection{Definition of the Brunerie number}
Brunerie's goal is now to analyse $\func{$\pi$}_3(\js)$. The first
result needed is the following:
\begin{defi}[Whitehead map]
  Given two pointed types $A$ and $B$, there is a map:
  \textnormal{
  \ExecuteMetaData[agda/latex/background.tex]{whitehead}
  }
\end{defi}
For our purposes, we only need the case when $A = B = \sphere{1}$
(although all of the following results appear in full generality in
Brunerie's thesis). We get a composite map:
\[
  \func{e} : \sphere{3}
  \xrightarrow{\,\func{$\simeq$}\,}
  {\JJoin{\sphere{1}}{\sphere{1}}}
  \xrightarrow{\,\func{W}\,}
  \sphere{2} \,\func{$\vee$}\, \sphere{2}
\]
This induces, via pre-composition, a \emph{Whitehead product}:
\[
  \func{$\pi$}_2(\sphere{2}) \,\func{$\times$}\, \func{$\pi$}_2(\sphere{2}) \xrightarrow{\,[-,-]\,} \func{$\pi$}_3(\sphere{2})
  \]
  by
  \[
[\trunc{f},\trunc{g}] := \trunc{\func{$\nabla$}\, \func{$\circ$}\, (f \,\func{$\vee$}\,g) \, \func{$\circ$}\, \func{e}}
  \]
Recall that we denote by $i_2$ the generator of
$\func{$\pi$}_2(\sphere{2})$.  Brunerie shows, in particular, the
following about its relation to the Whitehead product (see~\cite[Proposition 3.4.4.]{Brunerie16} for the full statement).
\begin{thm}\label{thm:main}
  The kernel of
  the suspension map $\func{$\sigma$}_*:\func{$\pi$}_3(\sphere{2}) \to
  \func{$\pi$}_4(\sphere{3})$ is generated by $[i_2,i_2]$.
\end{thm}
The key technical component in the proof is the \emph{Blakers-Massey
  Theorem}, first formalised in HoTT
by~Favonia,~Finster,~Licata~\&~Lumsdaine~in~\cite{FavoniaFinster+16}:
\begin{thm}[Blakers-Massey]\label{thm:bm}
  Consider the diagram
\[\begin{tikzcd}[ampersand replacement=\&]
	A \\
	\& P \& C \\
	\& B \& {\func{\textnormal{Pushout}}\,f\,g}
	\arrow[from=2-2, to=3-2]
	\arrow[from=2-2, to=2-3]
	\arrow["{\inl{}}"', from=3-2, to=3-3]
	\arrow["{\inr{}}", from=2-3, to=3-3]
	\arrow["f"', curve={height=12pt}, from=1-1, to=3-2]
	\arrow["g", curve={height=-12pt}, from=1-1, to=2-3]
	\arrow["{f\sqcup g}"{description}, from=1-1, to=2-2]
	\arrow["\lrcorner"{anchor=center, pos=0.125}, draw=none, from=2-2, to=3-3]
        \arrow["\lrcorner"{anchor=center, pos=0.125, rotate=180}, draw=none, from=3-3, to=2-2]
      \end{tikzcd}\]
    where $P$ is the pullback along $\inl{}$ and $\con{inr}$,
    i.e.\ $P = \Sigma_{(b,c) : B
      \,\func{$\times$}\,C}(\tyPath{\inl{b}}{\inr{c}})$, and $f \sqcup
    g$ is defined by
    $$(f \sqcup g)\,a = (f \,a \,\con{,}\, g\,a \,\con{,}\,\push{a})$$
    If $f$ and $g$ are $n$- respectively
    $m$-connected, then $f \sqcup g$ is $(n+m)$-connected.
\end{thm}
\autoref{thm:main} is proved by considering the following diagram
\[\begin{tikzcd}[ampersand replacement=\&]
	\sphere{3} \\
	\& P \& \sphere{2} \\
	\& \One \& {J_2\sphere{2}}
	\arrow[from=2-2, to=3-2]
	\arrow[from=2-2, to=2-3]
	\arrow[from=3-2, to=3-3]
	\arrow[from=2-3, to=3-3]
	\arrow[curve={height=12pt}, from=1-1, to=3-2]
	\arrow["{\Nabla \circ \,\func{W}}"{pos = 0.38 , rotate = -18}, curve={height=-12pt}, from=1-1, to=2-3]
	\arrow[from=1-1, to=2-2]
	\arrow["\lrcorner"{anchor=center, pos=0.125}, draw=none, from=2-2, to=3-3]
        \arrow["\lrcorner"{anchor=center, pos=0.125, rotate=180}, draw=none, from=3-3, to=2-2]
      \end{tikzcd}\]

Verifying that the outer square is a pushout square is technical and
we refer to Brunerie's proof for the details. Above, $P$ is simply the fibre of
$\inr{}:\sphere{2} \to \js$. The leftmost map is
$2$-connected since $\sphere{3}$ is $2$-connected
and the top map is
$0$-connected since $\sphere{3}$ and $\sphere{2}$ are both
$1$-connected. Consequently, by \autoref{thm:bm}, we get that the map
$\sphere{3}\to P$ is $2$-connected and thus induces a surjection after application of
$\func{$\pi$}_3$. This gives the diagram:
\[\begin{tikzcd}[ampersand replacement=\&]
	{\func{$\pi$}_3(P)} \& {\func{$\pi$}_3(\sphere{2})} \& {\func{$\pi$}_3(\js)}\\
	{\func{$\pi$}_3(\sphere{3})} \&\& \func{$\pi$}_4(\sphere{3})
	\arrow[from=1-1, to=1-2]
	\arrow[from=1-2, to=1-3]
        \arrow["\func{$\cong$}",from=1-3, to=2-3]
        \arrow["\func{$\sigma$}_*", swap, from=1-2, to=2-3]
	\arrow[from=2-1, to=1-2,dashed]
	\arrow[two heads, from=2-1, to=1-1]
\end{tikzcd}\]
where the sequence on the top comes from the long
exact sequence of homotopy groups associated to $P$. The dashed map
sends the generator $i_3:\func{$\pi$}_3(\sphere{3})$ to $[i_2,i_2] :
\func{$\pi$}_3(\sphere{2})$ by definition.

\autoref{thm:main} motivates the following definition. Recall that we
denote by $\func{$\psi$}$ the isomorphism
$\tyIso{\func{$\pi$}_3(\sphere{2})}{\bZ}$.
\begin{defi}[Brunerie number]
  We define the Brunerie number $\brunerie:\bZ$ by
  $\brunerie = \func{$\psi$}[i_2,i_2]$.
\end{defi}

We may now prove the main result of \cite[Chapter~3]{Brunerie16}.
\begin{cor}\label{cor:main}
  $\tyIso{\func{$\pi$}_4(\sphere{3})}{\bZ/\brunerie\bZ}$.
\end{cor}
\begin{proof}
  We have a homomorphism $\func{$\sigma$}_*
  \,\func{$\circ$}\,\func{$\psi$}^{-1} : \bZ \to
  \func{$\pi$}_4(\sphere{3})$. This composition is surjective since
  $\func{$\psi$}$ is an isomorphism and $\func{$\sigma$}_* :
  \func{$\pi$}_3(\sphere{2}) \to \func{$\pi$}_4(\sphere{3})$ is
  surjective by \autoref{thm:Freudenthal}. Since, by
  \autoref{thm:main}, the kernel of $\func{$\sigma$}_*$ is generated
  by $[i_2,i_2]$, the kernel of $\func{$\sigma$}_*
  \,\func{$\circ$}\,\func{$\psi$}^{-1}$ is generated by
  $\func{$\psi$}[i_2,i_2]$, i.e. by $\func{$\beta$}$. The statement then follows from the first isomorphism theorem.
\end{proof}

\subsection{Formalisation of the definition of the Brunerie number}

The formalisation of this part was straightforward. Arguably the most
technical result, the Blakers-Massey theorem, was already available in
the library thanks to~Kang~\cite{BMFormalization}. Most of the
remaining results were essentially just diagram chases which, in a
proof assistant, can be somewhat technical. Most work went into
verifying that $\js$ is the cofibre of $\Nabla \circ \func{W}$, the
proof of which followed Brunerie's closely.

In this section we found the only obvious mistake in Brunerie's
thesis. On page 82, in his definition of the $\push{}$\!-case for
$\func{W}$, the path component in the middle was not inverted, making
the term ill-typed. Naturally, this was of no mathematical
significance and something Brunerie immediately would have noticed if
he would have attempted to provide a computer formalisation of this
construction.

%% file: sections/brunerieproof.tex
\section{Brunerie's proof that $\tyPath{\abs{\brunerie}}{2}$}
\label{sec:brunerieproof}

This section concerns the final three Chapters (4--6) of Brunerie's
thesis. The main goal here is proving that $\tyPath{\abs{\brunerie}}{2}$.

We will not discuss Chapter 4 in much detail. Chapter 4 is devoted to
smash products and, in particular, their symmetric monoidal
structure. Brunerie used this in subsequent chapters to define and
prove properties about the \emph{cup product}, a graded multiplicative
operation on cohomology groups which will be used to show that
$\tyPath{\abs{\brunerie}}{2}$. This chapter has turned out to be
incredibly difficult to formalise due to the large number of higher
coherences involved in the proofs \cite{brunerie18}. In fact, the
results of this chapter were proved in detail and fully formalised
only $8$ years after the publication of Brunerie's thesis by
Ljungström~\cite{axel-smash}. We remark that despite the fact that
these results have now been made available to us, they are not
needed. While, with these results, Brunerie's construction of the cup
product appears correct, his use of smash products still leads to some
rather cumbersome diagram chases (with many coherences which still
need verification.)

Luckily, it turns out that Chapter 4 can be avoided altogether and
that this in fact makes some difficult proofs later on very
direct. For this reason, the results in Chapter~4 were omitted
completely from our formalisation. The reason for this is that all
results regarding smash products in Brunerie's thesis concern, in some
way, pointed maps out of smash products. In this case, we may exploit
the adjunction of maps out of smash products and bi-pointed maps:
\[
  \tyEquiv{(A\,\func{$\wedge$}\, B \to_\star C)}{(A \to_\star (B \to_\star C))}
\]
Here, $B \to_\star C$ is taken to be pointed by the constant map. As
shown by Brunerie et al.~\cite{BLM22}, it is arguably easier to define
the cup product on the right-hand side of the adjunction, which
effectively means that we never have to work with smash products when
formalising cohomology theory. The usefulness of the approach by
Brunerie et al.~\cite{BLM22} is not only witnessed by our work---it
has been used by Lamiaux et al.~\cite{LLM23} and
Ljungström~\&~Mörtberg~\cite{LM24} in the development of cohomology
rings and is used to describe the cup product as an instance of the
delooping machinery introduced by Wärn~\cite[Section 4.3]{david23}. We
remark that the same techniques (although independent
from~\cite{BLM22}) can be found in the work of
Christensen~\&~Scoccola~\cite[Section 2.4]{christensen2020hurewicz}
where it is utilised in a discussion of the magma structure on loop
spaces.

\subsection{Cohomology and the Hopf invariant}

\cite[Chapter~5]{Brunerie16} introduces integral cohomology groups and
rings, and gives a construction of the Mayer-Vietoris sequence. In
more detail, Brunerie defines the integral Eilenberg-MacLane spaces by
$\EM{0} = \bZ $ and $\EM{n} = \truncT{n}{\sphere{n}}$ for $n \geq
1$. This allows for a definition of the (integral) cohomology of $X$:
\[
  \cohom{n}{X} = \truncT{0}{X \to \EM{n}}
\]
The fact that $\tyEquiv{\LoopSpace{\EM{n+1}}}{\EM{n}}$ follows by a proof
completely analogous to that of~\autoref{cor:stable}. Brunerie uses this
equivalence to carry over the (commutative) H-space structure on
$\LoopSpace{\EM{n+1}}$ to that of $\EM{n}$. This provides a notion of
addition $\plusk : \EM{n} \,\func{$\times$}\, \EM{n} \to
\EM{n}$ which lifts to $\cohom{n}{X}$ by post-composition, thereby endowing
$\cohom{n}{X}$ with a group structure.

Similarly, Brunerie gives a definition of a cup product
$\cupprodkop : \EM{n}\, \Smash \,\EM{m} \to \EM{n+m}$ which lifts to
the usual cup product $\cupprodop : \cohom{n}{X}
\,\func{$\times$}\, \cohom{m}{X} \to \cohom{n+m}{X}$. This is shown to
induce a graded commutative ring structure on $\cohom{*}{X}$ using
results from Chapter~4.

The synthetic construction of the Mayer-Vietoris sequence concerns the
long exact sequence
\[\begin{tikzcd}[ampersand replacement=\&,row sep = -2pt, column sep = 38pt]
	{\cohom{0}{D}} \& {\cohom{0}{B}\times \cohom{0}{C}}  \& {\cohom{0}{A}} \\
	\phantom{.} \& \phantom{.} \& \phantom{.} \\
	{\cohom{1}{D}} \& \dots \&
	\arrow[from=1-1, to=1-2]
	\arrow[from=1-2, to=1-3]
        \arrow[from=3-1, to=3-2]
	\arrow[from=2-2, to=3-1,shorten <=-7.5pt, curve={height=6pt}]
	\arrow[shorten <=4pt,shorten >=-5pt, no head, from=1-3, to=2-2,curve={height=-3pt}]
\end{tikzcd}\]
where $D$ denotes the pushout of a span
$B \xleftarrow{\,f\,} A \xrightarrow{\,g\,} C$. A direct application gives us,
for $n \geq 1$, that $\tyIso{\cohom{n}{\sphere{m}}}{\bZ}$ if $n = m$
and $\tyIso{\cohom{n}{\sphere{m}}}{\One}$ otherwise.  This gives, by
another application of the sequence, the following result:
\begin{lem}\label{lem:cohom-push}
  For any $f : \sphere{3} \to \sphere{2}$ we have
  \textnormal{
  \begin{align*}
    \cohom{n}{\cofib{f}} \,\,\func{$\cong$}\,\, \begin{cases} \bZ & n \in \{0,2,4\} \\ \One &
      \text{otherwise}
      \end{cases}
  \end{align*}}
\end{lem}
Let us briefly fix $f : \sphere{3} \to \sphere{2}$. Denote by
$\gamma_2$ and $\gamma_4$ the generators of $\cohom{2}{\cofib{f}}$ and
$\cohom{4}{\cofib{f}}$ respectively given by the image of $1:\bZ$
under the isomorphism in~\autoref{lem:cohom-push}. These generators
may be used to define an invariant on $\sphere{3} \to \sphere{2}$
called the \emph{Hopf invariant}. This is done as follows:
\begin{defi}[Hopf invariant]\label{def:HopfInvariant}
  The \emph{Hopf invariant} of $f$ is the unique integer
  $\HI{f} : \bZ$ such that
  $\tyPath{\cupprod{\gamma_2}{\gamma_2}}{\HI{f}}\,\func{$\cdot$}\, \gamma_4$.
\end{defi}
We remark that the above definition is given for the more general
class of maps $\sphere{2n-1}\to\sphere{n}$ in Brunerie's thesis.
For our purposes, the above special case
suffices.
In particular, we may see $\textnormal{\func{HI}}$ as a function
$\pi_3(\sphere{2}) \to \bZ$. The following turns out to be true:
\begin{prop}\label{prop:HI-homomorphism}
  $\textnormal{\func{HI}}$ is a homomorphism
  $\pi_3(\sphere{2}) \to \bZ$.
\end{prop}
\begin{proof}[Proof sketch]
  We first rephrase $f \,\func{+}\, g : \pi_3(\sphere{2})$ as a
  composition
  \[
    \sphere{3}
    \xrightarrow{}
    \sphere{3}\,\func{$\vee$}\,\sphere{3}
    \xrightarrow{\,f \,\func{$\vee$}\, g\,}
    \sphere{2}\,\func{$\vee$}\,\sphere{2}
    \xrightarrow{\func{$\nabla$}}
    \sphere{2}
  \]
  By analysing the cohomology of $\cofib{(\func{$\nabla$}\,\func{$\circ$}\,(f \,\func{$\vee$}\, g))}$ and
  the action on generators of the obvious maps from
  $\cofib{(\func{$\nabla$}\,\func{$\circ$}\,(f \,\func{$\vee$}\, g))}$, $\cofib{f}$ and $\cofib{g}$ into
  $\cofib{(f \,\func{+}\, g)}$, one arrives at the result with some
  elementary algebra.
\end{proof}

Finally, the Hopf invariant of our element of interest $[i_2,i_2]$ is
computed (up to a sign), using an argument similar to that of the
proof of~\autoref{prop:HI-homomorphism}.
\begin{prop}\label{prop:HI-whitehead}
  $\tyPath{\abs{\HI{[i_2,i_2]}}}{2}$
\end{prop}

We are now almost done: if there is an element $f : \pi_3(\sphere{2})$
such that $\tyPath{\HI{f}}{1}$, then $\HI{}$ is an isomorphism
$\tyIso{\pi_3(\sphere{2})}{\bZ}$. Since isomorphisms of this type are
unique up to a sign, \autoref{prop:HI-whitehead} tells us that also
for the standard isomorphism
$\func{$\psi$} : \tyIso{\pi_3(\sphere{2})}{\bZ}$, we must have
$\tyPath{\abs{\func{$\psi$}[i_2,i_2]}}{2}$, i.e.\
$\tyPath{\abs{\brunerie}}{2}$. Hence, we have so far shown the
following:
\begin{lem}\label{lem:almost}
  If $\tyPath{\HI{f}}{1}$ for some $f : \pi_3(\sphere{2})$, then $\tyPath{\abs{\brunerie}} {2}$.
\end{lem}
The final chapter of Brunerie's thesis is devoted to proving the
antecedent of \autoref{lem:almost}.

\subsection{Formalisation of cohomology and the Hopf invariant}

This section was largely covered by Brunerie, Ljungström and Mörtberg
in~\cite{BLM22} and thus also available in \agdaCubical. We briefly
summarise:
\begin{itemize}
\item $\func{$+_k$} : \EM{n} \times \EM{n} \to \EM{n}$ was defined
  explicitly using a direct construction of the \emph{Wedge
    Connectivity} for spheres---a generalisation
  of~\autoref{lem:wedge-conn}. This construction is of great convenience
  to our formalisation due to the fact that e.g.\
  $\star_{\EM{n}}\,\func{$+_k$}\,\trunc{x}\,\func{$\equiv$}\,\trunc{x}$
  holds definitionally. In fact, all of the basic laws governing
  $\func{$+_k$}$ are (trivially) provably path equal to $\refl$ at
  $\star_{\EM{n}}$, which simplifies a lot of path algebra.
  \item The cup product is defined
    via the following lift
    \[\begin{tikzcd}[ampersand replacement=\&]
	{\sphere{n}} \&\& {(\EM{m} \to_\star \EM{n+m})} \\
	{\EM{n}}
	\arrow["{\trunc{\!\con{\_\,}\!}}"', from=1-1, to=2-1]
	\arrow["\cupprodop"', dashed, from=2-1, to=1-3]
	\arrow[from=1-1, to=1-3]
      \end{tikzcd}\]
    for $n \geq 1$, where the top map may be thought of as being inductively defined
    via the equivalence
    \begin{align*}
      &\phantom{\func{$\simeq$}\,} (\sphere{n+1} \to_\star (\EM{m} \to_\star \EM{(n+1)+m})) \\
      &\func{$\simeq$}\, (\sphere{n} \to_\star (\EM{m} \to_\star \LoopSpace{\EM{(n+1)+m}})) \\
      &\func{$\simeq$}\, (\sphere{n} \to_\star (\EM{m} \to_\star \EM{n+m}))
    \end{align*}
    The lift exists because the type of pointed functions
    $\EM{m} \to_\star \EM{n+m}$ is an $n$-type. This construction
    gives an inductively defined cup product which is remarkably easy
    to work with, as showcased in \cite{LLM23} to compute cohomology
    rings of various classical spaces.
  \item The Mayer-Vietoris sequence was formalised by directly
    translating Brunerie's original proof.
\end{itemize}
Hence, what remained to be formalised in Chapter 5 was the Hopf
invariant, \autoref{prop:HI-homomorphism} and
\autoref{prop:HI-whitehead}. The formalisation of these propositions
was straightforward and we were able to translate Brunerie's proofs in
a direct manner. This is not surprising as the proofs are very
algebraic.

For simplicity, we only formalised these propositions as they stand
here and not their generalisations to higher spheres (i.e.\ as in
\cite[Proposition~5.4.3~\&~5.4.4]{Brunerie16}). We remark, however,
that the formalised proofs easily should be rephrasable for the
general Hopf invariant of maps $\sphere{2n-1} \to \sphere{n}$.

\subsection{The Gysin sequence}

This section corresponds to \cite[Chapter~6]{Brunerie16}. In order to
be able to apply \autoref{lem:almost}, this chapter is devoted to
proving that $\tyPath{\abs{\HI{\,\hopf}}}{1}$, where, recall,
$\hopf : \sphere{3} \to \sphere{2}$ is the Hopf map---the generator of
$\pi_3(\sphere{2})$ from \autoref{def:hopfmap}. This amounts to
analysing the cup product on the cohomology of $\cofib{\hopf}$. It is
well-known that $\cofib{\hopf}$ is a model of the complex
projective plane $\CP$ (see e.g.\ \cite[Example 4.45]{Hatcher2002}),
so let us simply write $\CP$ from now on. We hence have $\CP$ defined
as the following pushout:
\[\begin{tikzcd}[ampersand replacement=\&]
	{\sphere{3}} \& {\sphere{2}} \\
	\One \& \CP
	\arrow["\hopf", from=1-1, to=1-2]
	\arrow[from=1-1, to=2-1]
	\arrow[from=2-1, to=2-2]
	\arrow[from=1-2, to=2-2]
	\arrow["\lrcorner"{anchor=center, pos=0.125, rotate=180}, draw=none, from=2-2, to=1-1]
\end{tikzcd}\]
In order to show that $\tyPath{\abs{\HI{\,\hopf}}}{1}$, it suffices to
show that $\cupprod{-}{\gamma_2} : \cohom{2}{\CP} \to \cohom{4}{\CP}$ is an isomorphism
for $\gamma_2 : \cohom{2}{\CP}$ a generator.
Brunerie does this by constructing the Gysin sequence.
\begin{prop}[The Gysin sequence]\label{prop:Gysin}
  Let $B$ be a pointed and $0$-connected type and $P : B \to
  \Type$
  be a fibration with $P\,\star_B \,\func{$\simeq$}\, \sphere{n-1}$.
  Let $E = \Sigma_{b : B} (P\,b)$ be the total space of $P$. If there is a family of maps
  $c : (b: B) \to (\Susp{(P\,b)} \to_\star \EM{n})$ with $c_{\star_B}$ a
  generator of $\cohom{n}{\sphere{n}}$, then there is an element
  $e_n : \cohom{n}{B}$ and a long exact sequence
  \[
\begin{tikzcd}[column sep=38pt,row sep = 6pt]
  & \dots\arrow[r] &
  \cohom{i-1}{B} \arrow[dll,out=280,in=30,looseness=0.2] \\
  {\cohom{i-1}{E}} \arrow[r] &
  {\cohom{i-n}{B}} \arrow["{\cupprod{-}{e_{n}}}",r] &
  {\cohom{i}{B}} \arrow[dll,out=280,in=30,looseness=0.2] \\
  {\cohom{i}{E}} \arrow[r] & \dots
\end{tikzcd}
\]
Moreover, $c$ (and also $e_n$) exists when $B$ is $1$-connected.
\end{prop}

In order to make use of this, we need the following result.
\begin{prop}\label{prop:CP-fib}
  There is a fibration $P : \CP \to \Type$ with $P\,\star_{\CP}
  \func{$\simeq$} \,\sphere{1}$ and total space $\sphere{5}$.
\end{prop}

\autoref{prop:CP-fib} is a special case of the following result.
\begin{prop}[Iterated Hopf construction]\label{prop:itHopf}
  Given an associative H-space $A$, let $h_A : \JJoin{A}{A} \to \Susp{A}$
  denote the associated Hopf map. There is a fibration
  \textnormal{$\cofib{h_A} \to \Type$ with fibre $A$ and
    total space $\JJoin{A\,}{\JJoin{A\,}{A}}$.  }
\end{prop}
We consider the particular case when $A = \sphere{1}$ in
\autoref{prop:itHopf}. In this case, the map $h_{\sphere{1}} : \JoinS \to
\sphere{2}$ corresponds to the usual Hopf map under the equivalence
$\tyEquiv{\JoinS}{\sphere{3}}$ and hence
$\tyEquiv{\cofib{h_{\sphere{1}}}}{\CP}$. The total space of this is
$\JJoin{\JoinS}{\sphere{1}}$ which is equivalent to $\sphere{5}$ by
\autoref{prop:spherejoin} and thus we have proved \autoref{prop:CP-fib}. The
associated Gysin sequence gives us the main result of this section:
\begin{prop}[Hopf invariant of the Hopf map]
  \label{prop:HI-1}
  $\tyPath{\abs{\HI{\hopf}}}{1}$
\end{prop}
\begin{proof}
  Since $\CP$ is $1$-connected, \autoref{prop:Gysin} combined with
  \autoref{prop:CP-fib} gives us an element $e_2 : \cohom{2}{\CP}$ and
  a sequence
  \[
    \cohom{i-1}{\sphere{5}} \to \cohom{i-2}{\CP} \xrightarrow{\,\cupprod{-}{e_2}\,} \cohom{i}{\CP} \to \cohom{i}{\sphere{5}}
  \]
  When $1 \leq i \leq 4$, $\cohom{i}{\sphere{5}}$ vanishes. Setting $i
  = 2$, we get that $e_2$ must be a generator of $\cohom{2}{\CP}$, and
  thus equal to the generator $\gamma_2 : \cohom{2}{\CP}$ up to a
  sign. Setting $i = 4$, we get that $\cupprod{-}{e_2}$ must be an
  isomorphism of groups $\tyIso{\cohom{2}{\CP}}{\cohom{4}{\CP}}$ and hence $\cupprod{e_2}{e_2}$ is a
  generator. Consequently, so is $\cupprod{\gamma_2}{\gamma_2}$, and thus $\tyPath{\abs{\HI{\,\hopf}}}{1}$.
\end{proof}

\autoref{prop:HI-1} combined with \autoref{lem:almost} gives the desired
path: $\tyPath{\abs{\brunerie}}{2}$. This completes Brunerie's proof
and \autoref{cor:main} gives us the main result:
\begin{thm}
  $\tyIso{\func{$\pi$}_4(\sphere{3})}{\bZ/2\bZ}$
\end{thm}

\subsection{Formalisation of the Gysin sequence}

Formalising the results from Chapter~6 was more challenging, but was
greatly aided by the alternative construction of the cup product
discussed above. The first technical lemma, which is
crucial for the construction of the Gysin sequence is:
\begin{lem}\label{lem:ap-cup}
  Given $x : \EM{n}$ and $y : \EM{m}$, we have\textnormal{
    \[\ap{(\lambda\,a \to \cupprodk{a}{y})}{(\func{$\sigma$}_{n}\,x)} \,\func{$\equiv$}\, \func{$\sigma$}_{n+m}(\cupprodk{x}{y})\]
  }
\end{lem}
In Brunerie's thesis, this lemma relies on a result which in turn
requires the symmetric monoidal structure of the smash product (in
particular, it uses the \emph{pentagon identity}). With the
alternative construction of the cup product, however, this result
follows immediately from the definition of the cup product.

\autoref{lem:ap-cup} is used to show that the map
\begin{align*}
  &g^{i} : \EM{i} \to (\sphere{n} \to_\star \EM{i+n}) \\
  &g^{i}\,x = \lambda y \to \cupprodk{x}{\iota{y}}
\end{align*}
is an equivalence, which is crucially used in the construction of the
Gysin sequence. Above, $\iota : \sphere{n} \to \func{K}_n$ is a
generator of $\cohom{n}{\sphere{n}}$. For reference, $g^i$ is the map
$g^{i}_{\star_{B}}$ in the proof of
\cite[Proposition~6.1.2]{Brunerie16}.
Brunerie's proof proceeds by induction
on $i$: the fact that $g^0$ is an equivalence is easy; for the
inductive step, it suffices to show that
$\LoopSpace{g^{i+1}} : {\LoopSpace {\EM{i+1}}} \to
{\LoopSpace{(\sphere{n} \to_\star \EM{(i+1)+n})}}$ is an equivalence
for reasons of connectedness. This is done by showing that the
following diagram commutes
\[\begin{tikzcd}[ampersand replacement=\&]
	{\LoopSpace{\EM{i+1}}} \&\& {\LoopSpace{(\sphere{n} \to_\star \EM{(i+1)+n})}} \\
	{\EM{i}} \&\& {(\sphere{n} \to_\star \EM{i+n})}
	\arrow["{\LoopSpace{g^{i+1}}}", from=1-1, to=1-3]
	\arrow["{{g^{i}}}", from=2-1, to=2-3,swap]
	\arrow["\func{$\simeq$}"', from=1-3, to=2-3,swap]
	\arrow["\func{$\simeq$}"', from=1-1, to=2-1]
\end{tikzcd}\]
hence getting that $\LoopSpace{g^{i+1}}$ is an equivalence from the
induction hypothesis.

While the general idea of Brunerie's proof of this statement is
correct, it was difficult to formalise directly. The primary reason
for this is that Brunerie does not pay much attention to the fact that
the objects of interest are not just functions, but \emph{pointed}
functions.  In particular, his argument for the commutativity of the
diagram above treats $g^i\,x$ as a plain function rather than a
pointed function.  Fortunately for us, the whole proof is very direct
with the alternative definition of the cup product. Formalising
Brunerie's proof with pointedness of functions respected would have
been hard, especially without machinery external to \cite{Brunerie16}
(e.g.\ \cite[Lemma 14.]{BLM22}).

After these subtleties were dealt with, the formalisation of the Gysin
sequence could proceed following Brunerie's proof closely. In our
initial formalisation, we made a slight adjustment to the indexing of
the Gysin sequence. This removed some bureaucracy but happened at the
cost of generality.\footnote{A more general form of the Gysin sequence
  using Brunerie's indexing has later been added to \agdaCubical.}
This made verifying that~\autoref{prop:HI-1} slightly less direct,
because we no longer had access to the case
\[
  \cohom{1}{\sphere{5}} \longrightarrow
  \cohom{0}{\CP} \xrightarrow{\,\cupprod{-}{e_2}\,}
  \cohom{2}{\CP} \longrightarrow
  \cohom{2}{\sphere{5}}
\]
which is used by Brunerie to show that the element
$e_2 : \cohom{2}{\CP}$, for which
$\cupprod{-}{e_2} : \cohom{2}{\CP} \to \cohom{4}{\CP}$ is an
isomorphism, is indeed a generator. However, in practice, this is not
a big problem. In fact, it provides a nice example of a proof by
computation. It is very direct to manually show that the map
$i : \CP \to \EM{2}$ induced by $i(\inl{x}) = \trunc{x}$ is equal to
the underlying map of $e_2$. The fact that $i$ generates
$\cohom{2}{\CP}$ can then be verified by computation: applying the
isomorphism $\tyIso{\cohom{2}{\CP}}{\bZ}$ to $\trunc{i}$ returns $1$
by normalisation in \CubicalAgda. We stress, for those skeptical of
this method, that it also is very direct to provide a ``manual''
formalisation of this fact.

The final step of the formalisation was~\autoref{prop:itHopf}, i.e.\ the
iterated Hopf construction.  Although technical, the formalisation
could be carried out following Brunerie closely.

%% file: sections/computation.tex
\section{The simplified proof and normalisation of a Brunerie number}
\label{sec:computation}

It turns out that not only Chapter 4, but also Chapters 5--6 can be
avoided. As conjectured by Brunerie, it would be possible to do this
by simply normalising the Brunerie number. While we still cannot
normalise his original definition of it, we can at least provide a
computation of a substantially simplified Brunerie number. This is
defined via a more tractable description of the isomorphism
$\tyIso{\func{$\pi$}_3(\sphere{2})}{\bZ}$ as a composition of simpler
isomorphisms, relying on an alternative definition of $\func{$\pi$}_3$
in terms of $\JoinS$. The idea is then to trace
$[i_2,i_2] : \func{$\pi$}_3(\sphere{2})$ step by step through these
isomorphisms. This gives a sequence of new Brunerie numbers and one of
these normalises to $- 2$ in \CubicalAgda in a matter of seconds.

The trick to give a more tractable definition of
$\tyIso{\func{$\pi$}_3(\sphere{2})}{\bZ}$ is to redefine the third
homotopy group of a type $A$ as $\func{$\pi$}_3^*(A) =
\truncT{0}{\JoinS \to_\star A}$. This reformulation of
$\func{$\pi$}_3$ can be given an explicit group structure, such that
pre-composition by $\tyEquiv{\JoinS}{\sphere{3}}$ induces an
isomorphism $\tyIso{\func{$\pi$}_3(A)}{\func{$\pi$}_3^*(A)}$. Let us
first set up machinery we need (and a bit more).

\subsection{Interlude: joins and smash products of spheres}
\label{interlude}
\newtheorem*{remark}{Remark}

We have seen that the equivalence
$\sphere{3}\,\func{≃}\,\JJoin{\sphere{1}}{\sphere{1}}$ played a
crucial role in Brunerie's original proof. What is less clear,
however, is what this equivalence actually looks like. It turns out
that it is closely related to the multiplication
$\sphere{1}\,\func{$\wedge$}\,\sphere{1} \to \sphere{2}$ and, as such,
has a rather direct and algebraic description. Let us therefore
briefly study this multiplication and describe its relation to the
decomposition of spheres into joins. Although we only need
low-dimensional special cases of these facts, we take the opportunity
to tell the general story.

\begin{remark}
  In this subsection only, we will use the definition
  \textnormal{$\sphere{n}:= \SuspN{n}{\func{Bool}}$}. In particular,
  this means that we redefine $\sphere{1} := \Susp{\func{Bool}}$
  instead of using the $\base/\Loop$ construction. This is only done
  for ease of presentation and is not used in the formalisation.
\end{remark}

Let us use the following (explicit) definition of the smash product
$A\,\func{$\wedge$}\,B$.
\ExecuteMetaData[agda/latex/background.tex]{smash}
This construction is well-known and can easily be seen to be (bi-)functorial (see e.g.~\cite[Definition
  6]{axel-smash}),
i.e.\ given pointed maps $f : A \to_\star C$ and $g : B \to_\star D$,
there is a map $f \,\Smash\, g : A \,\Smash\,B\to C\,\Smash D$ with
$(f \,\Smash\, g)\smashelem{x}{y} := \smashelem{f\,x}{g\,y}$.

The first goal is to define a multiplication $\sphere{n}\, \Smash \,
\sphere{m} \to \sphere{n+m}$. To facilitate future proofs, we first
introduce the following construction which lifts maps $A
\,\func{$\times$}\, B \to C$ to maps $(\Susp{A}) \,\func{$\times$}\,B
\to \Susp{B}$:
\ExecuteMetaData[agda/latex/background.tex]{lift-left}
The function $\liftL{f}$ is pointed in the left-argument by
construction. It is pointed also in the right-argument if this is the
case for $f$. Hence, given any function $g : A \,\Smash\,B \to C$, we
also get (with some abuse of notation) $\liftL{g} : (\Susp{A})
\,\Smash\,B \to \Susp{C}$.
\begin{lem}\label{lem:susp-smash}
  For any pointed types $A$ and $B$, the map $\liftL{\,\idfun}_{A
    \,\Smash\,B} : (\Susp{A}) \,\Smash\,B \to \Susp{(A \,\Smash\,B)}$ is
  an equivalence.
\end{lem}
\begin{proof}
  The inverse of $\liftL{\,\idfun}_{A \,\Smash\,B}$ is induced by the
  map $A \,\func{$\times$}\, B \to \LoopSpace{((\Susp{A}) \,\Smash\,B})$ defined by
  mapping $(a,b) : A \,\func{$\times$}\, B$ to the composite loop given by:
\[
 \begin{tikzcd}[ampersand replacement=\&]
	{\smashpt} \& {\langle \north,b\rangle} \& \& {\smashelem{\north}{b}} \&{\smashpt}
	\arrow["\pushr^{\func{$-1$}}", from=1-1, to=1-2]
	\arrow["{\,\ap{\smashelem{\func{$-$}}{b}}{(\func{$\sigma$}\,a)}\,}", from=1-2, to=1-4]
	\arrow["\pushr", from=1-4, to=1-5]
 \end{tikzcd}
 \]
The fact that these maps cancel follows by some technical but
elementary path algebra. For the details, we refer to the
formalisation.
\end{proof}
\begin{lem}
  If $f : A \,\Smash\,B \to C$ is an equivalence, then so is
  $\liftL{f} : \Susp{A}\,\Smash\,B \to \Susp\,A$
\end{lem}
\begin{proof}
  Using equivalence induction~(see
  e.g.~\cite[Corollary~5.8.5]{HoTT13}), it is enough to prove the
  lemma for $C := A \,\Smash\,B$ and $f := \idfun_{A\,\Smash\,B}$.
  In this case, the statement is precisely that of~\autoref{lem:susp-smash}.
\end{proof}
These lemmas allow us to define an equivalence $\func{$\wedge$}_{n,m}
: \sphere{n} \,\Smash\, \sphere{m} \simeq \sphere{n+m}$ (we borrow
this notation from~\cite[Proposition 4.2.2]{Brunerie16}). We will
write $\func{$\smile$} : \sphere{n} \, \func{$\times$}\, \sphere{m}
\to \sphere{n+m}$ for the underlying function, i.e. $\cupprod{x}{y} :=
\func{$\wedge$}_{n,m}\smashelem{x}{y}$.  The name $\func{$\smile$}$ is
suggestive: modulo the quotient maps $\sphere{\bullet} \to
\func{K}_{\bullet}$, it is precisely the cup product; this justifies
overloading the symbol. We define it by induction on $n$. In the case
$n = 0$, we define it on canonical points $\smashelem{x}{y}$ by case
distinction on $x$:
\begin{align*}
  \cupprod{\falsebool}{y} &= y\\
  \cupprod{\truebool}{y} &= \star_{\sphere{m}}
\end{align*}
This map induces a map on the full smash product
$\sphere{0}\,\Smash\,\sphere{m}$ \cite[Section~4.1]{Brunerie16}. In
fact, it is an equivalence, and thereby $\func{$\wedge$}_{0,m}$ is
defined. For $n>0$ we use the fact that $\sphere{n} :=
\Susp{\sphere{n-1}}$ and simply define define $\func{$\wedge$}_{n,m}
:= \liftL{\func{$\wedge$}_{n-1,m}}$. By \autoref{lem:susp-smash}, this
is an equivalence.

Let us try to transfer this construction from smash products to
joins. To begin with, consider the following map, defined for any two
pointed types $A$ and $B$:
\ExecuteMetaData[agda/latex/background.tex]{pinch}
\begin{prop}\label{prop:pinch-equiv}
For any two pointed types $A$ and $B$, the map $\pinch$ is an equivalence.
\end{prop}
\begin{proof}
While we have, in our formalisation, explicitly constructed an inverse
of $\pinch$ and proved directly that the two maps cancel, a recent
(independent) result by Cagne et al.~\cite{cagne2024symmetries} allows
us to give a more principled proof. We proceed by noting that for any
pointed type $C$, we have
\[
  (\Susp{(A \,\Smash\,B)} \to_\star C)
  \,\func{$\simeq$}\, ({A \,\Smash\,B} \to_\star \LoopSpace{C})
  \,\func{$\simeq$}\, (A \to_\star (B\to_\star \LoopSpace{C}))
  \,\func{$\simeq$}\, (\JJoin{A\,}{B} \to_\star C)
\]
where the first equivalence comes from the adjunction between
\func{Susp} and $\func{$\Omega$}$ and the second from the adjunction
between smash products and doubly pointed maps. The third equivalence
is \cite[Lemma 6.1]{cagne2024symmetries}. This shows that $\Susp{(A
  \,\Smash\,B)}$ and $\JJoin{A\,}{B}$ have the same elimination
principle, which implies the desired statement.
\end{proof}

Let $\func{F}_{n,m} : \JJoin{\sphere{n}}{\sphere{m}} \to
\sphere{n+m+1}$ denote the following composition:
\[
\JJoin{\sphere{n}}{\sphere{m}}
\xrightarrow{\,\pinch\,}
\Susp{(\sphere{n}\,\Smash\,\sphere{m})}
\xrightarrow{\,\Susp{(\Smash_{n,m})}\,}
\Susp{\sphere{n+m}} =: \sphere{n+m+1}
\]
Unfolding the definition of $\func{F}_{n,m}$ we see that it has an
incredibly compact description:
\ExecuteMetaData[agda/latex/background.tex]{F}
Since $\func{F}_{n,m}$ is a composition of two equivalences, we
immediately arrive at the following result.
\begin{prop}\label{prop:F-equiv}
\textnormal{$\func{F}_{n,m}$} is an equivalence
\end{prop}

We have already seen that the special case
$\JJoin{\sphere{1}}{\sphere{1}} \,\func{$\simeq$}\,\sphere{3}$ plays
an important role in Brunerie's proof. Now that we have moved from the
rather opaque definition of this equivalence presented in Brunerie's
thesis to the definition in terms of the very explicit function
$\func{F}_{n,m}$, we can hope to better understand its role in the
definition of the Brunerie number. Since the only non-trivial
component of the construction of $\func{F}_{n,m}$ is $\cupprodop$, we
may hope that it inherits some of its properties. We study these now.

Let us first analyse the interaction of $\cupprodop$ with
inversion. Recall, given a pointed type $A$ we can define inversion on
$\Susp{A}$ by:
\ExecuteMetaData[agda/latex/background.tex]{susp-inversion}
We get sphere inversion by letting $- : \sphere{n} \to \sphere{n}$ be
boolean negation when $n := 0$ and the suspension inversion defined
above when $n>1$.\footnote{If we prefer to use the
  $\base/\Loop$-construction of $\sphere{1}$, we may define the
  inversion map simply by sending $\Loop$ to $\Loop^{\func{$-1$}}$.}
\begin{prop}\label{prop:cup-comm}
The multiplication $\cupprodop$ is graded-commutative, i.e.\ for $x :
\sphere{n}$ and $y : \sphere{m}$, we have $\cupprod{x}{y}
\,\func{$\equiv$}\,\func{$-$}^{nm} (\cupprod{y}{x})$.
\end{prop}
For the proof, we refer to \cite[Proposition 18]{BLM22} which is the
corresponding statement for the cup product on $\func{K}_n :=
\truncT{n}{\sphere{n}}$ and whose proof directly applies also in our
setting. Associativity follows, just like in the proof of
\cite[Proposition 17]{BLM22}, by sphere induction:
\begin{prop}\label{prop:assoc}
The multiplication $\cupprodop$ is associative.
\end{prop}
\begin{proof}
Let $x : \sphere{n}$, $y : \sphere{m}$ and $z:\sphere{k}$. We show
that $\cupprod{x}{(\cupprod{y}{z})}
\,\func{$\equiv$}\,\cupprod{(\cupprod{x}{y})}{z}$ by induction on $n$
and $x$. When $n = 0$, the two equalities
\begin{align*}
\cupprod{\falsebool}{(\cupprod{y}{z})} \,\func{$\equiv$}\,\cupprod{(\cupprod{\falsebool}{y})}{z}\\
\cupprod{\truebool}{(\cupprod{y}{z})} \,\func{$\equiv$}\,\cupprod{(\cupprod{\truebool}{y})}{z}
\end{align*}
hold definitionally. For the inductive step, we use suspension
elimination on $x : \sphere{n+1}$. The two equalities
\begin{align*}
\cupprod{\north}{(\cupprod{y}{z})} \,\func{$\equiv$}\,\cupprod{(\cupprod{\north}{y})}{z}\\
\cupprod{\south}{(\cupprod{y}{z})} \,\func{$\equiv$}\,\cupprod{(\cupprod{\south}{y})}{z}
\end{align*}
also hold definitionally. So, by inspection of the definition of
$\cupprodop$, we need to show that
\begin{align*}
\func{$\sigma$}(\cupprod{x}{y})\,\func{$\equiv$}\,\ap{(\cupprod{-}{z})}{(\func{$\sigma$}{(\cupprod{x}{y})})}
\end{align*}
Using the action of \func{cong} on path composition, we can unfold the
right-hand side as follows:
\begin{align*}
\ap{(\cupprod{-}{z})}{(\func{$\sigma$}{(\cupprod{x}{y})})} \,&\func{$\equiv$}\, \ap{(\cupprod{-}{z})}{(\merid{(\cupprod{x}{y})})} \,\func{$\cdot$}\, \ap{(\cupprod{-}{z})}{(\merid{\north})^{\func{$-1$}}}  \\
&{:=}\, \func{$\sigma$}(\cupprod{x}{y}) \,\func{$\cdot$} \, \func{$\sigma$}(\cupprod{\north}{y})^{\func{$-1$}}\\
&\func{$\equiv$}\, \func{$\sigma$}(\cupprod{x}{y}) \qedhere
\end{align*}
\end{proof}

\subsection{Homotopy groups in terms of joins}
\label{htpyintermsofjoins}

As we have seen in Brunerie's construction of the Hopf map, it is
often easier to describe maps of type $\JJoin{\sphere{n}}{\sphere{m}}
\to A$ than those of type $\sphere{n+m+1} \to A$. However, the
definition of homotopy groups we have relied on so far uses the latter
type. This forces us to translate back and forth whenever we want to
use the definition in terms of joins. The key strategy behind our new
calculation of the Brunerie number is to rephrase homotopy groups in
terms of maps out of joins of spheres.
\begin{defi}
Given a pointed type $A$, we define
$\func{$\pi$}_{n+m+1}^{*}(A):=\truncT{0}{\JJoin{\sphere{n}}{\sphere{m}}
  \to_\star A}$.
\end{defi}
Clearly, this is equivalent to the usual definition of
$\func{$\pi$}_{n+m+1}(A)$ via pre-composition by
$\func{F}_{n,m}$. However, $\func{$\pi$}_{n+m+1}^{*}(A)$ can be
endowed with an explicit group structure which $\func{F}_{n,m}$ turns
out to respect. In order to construct the group structure on
$\func{$\pi$}^*_{n+m+1}(A)$, let us construct a map $\func{$\ell$} : A
\,\func{$\times$}\,B \to \LoopSpace{(\JJoin{A\,}{B})}$ for all pointed
types $A$ and $B$. Recall, we take $\JJoin{A\,}{B}$ to be pointed by
$\inl{\star_A}$. We define \func{$\ell$} by:
\[ \func{$\ell$}(a\,\con{,}\,b) := \push{(\star_A\,\con{,}\,\star_B)}\,\func{$\cdot$}\, \push{(a\,\con{,}\,\star_B)}^{\func{$-1$}} \,\func{$\cdot$}\, \push{(a\,\con{,}\,b)} \,\func{$\cdot$}\, \push{(\star_A\,\con{,}\,b)}^{\func{$-1$}} \]
Note that $\func{$\ell$}$ is pointed in both arguments. Let us also
define explicitly (once and for all) a pointed version of
$\func{cong}$ taking a pointed functions $f : A \to_\star B$ to a
pointed function $\func{$\text{cong}_\star$}\,f : \LoopSpace{A}
\to_\star \LoopSpace{B}$. We define it by
\[ \congpt{f}{p} := \star_f^{\func{$-1$}}\,\func{$\cdot$}\, \ap{f}{p} \,\func{$\cdot$}\,\star_f\]
where, recall, $\star_f : f\,\star_A \,\func{$\equiv$}\, \star_B$.
In other words, \func{$\text{cong}_\star$} is the functorial action
of $\func{$\Omega$}$.

We can now add two functions $f$ and $g$ of type $\JJoin{A\,}{B}
\to_\star C$ by
\ExecuteMetaData[agda/latex/background.tex]{joinMult}
We take this function to be pointed by $\refl$. Note that, since
\func{$\ell$} is pointed in both arguments, both
$\ap{(f\,\func{$+^{*}$}\,g)}{(\push{(a\,\con{,}\,\star_B)})}$ and
$\ap{(f\,\func{$+^{*}$}\,g)}{(\push{(\star_A\,\con{,}\,b)})}$
vanish. Let us compare this with the addition on the usual definition
homotopy groups. In general, we may add any two functions $f$ and $g$
of type $\Susp{A} \to_\star B$ by
\ExecuteMetaData[agda/latex/background.tex]{suspMult}
This is precisely the construction used to define the group structure
on $\func{$\pi$}_{n}$ whenever $n > 0$. Note that, by construction, we
have
$\ap{(f\,\func{$+^{\func{\textnormal{Susp}}}$}\,g)}{(\merid{\star_A})}
\,\func{$\equiv$}\, \refl$.

\begin{prop}\label{prop:+-pres}
Given $f, g : \sphere{n+m+1} \to_\star A$, we have
\textnormal{
  \[(f\,\func{$+^{\func{\textnormal{Susp}}}$}\,g)
  \,\func{$\circ$}\, \func{F}_{n,m} \,\func{$\equiv$}\,
  (f\,\func{$\circ$}\,\func{F}_{n,m}) \func{$+^{*}$}
  (g\,\func{$\circ$}\,\func{F}_{n,m})
  \]}
\end{prop}
\begin{proof}
The two functions agree on $\con{inl}$ and $\con{inr}$ by $\refl$. Let us consider the action on $\con{push}(x\,\con{,}\,y)$. We have\textnormal{
\begin{align*}
\ap{((f\,\func{$+^{\func{\textnormal{Susp}}}$}\,g)\,\func{$\circ$}\,\func{{F}}_{n,m})}{(\push{(x\,\con{,}\,y)})} \,&{:=} \,\ap{(f\,\func{$+^{\func{\textnormal{Susp}}}$}\,g)}{(\func{$\sigma$}(\cupprod{x}{y}))}\\
&\func{$\equiv$}\,\ap{(f\,\func{$+^{\func{\textnormal{Susp}}}$}\,g)}{(\merid{(\cupprod{x}{y})})} \\
&\phantom{.}\func{$\cdot$}\,\, \ap{(f\,\func{$+^{\func{\textnormal{Susp}}}$}\,g)}{(\merid{\north})}^{\func{$-1$}} \\
&\func{$\equiv$}\,\ap{(f\,\func{$+^{\func{\textnormal{Susp}}}$}\,g)}{(\merid{(\cupprod{x}{y})})}
\end{align*}}
which, by definition, unfolds to
\begin{align}
\label{eq:+susp}
\congpt{f}{(\func{$\sigma$}(\cupprod{x}{y}))}\,\func{$\cdot$}\, \congpt{g}{(\func{$\sigma$}(\cupprod{x}{y}))}
\end{align}
On the other hand, $\ap{((f\,\func{$\circ$}\,\func{F}_{n,m}) \func{$+^{*}$} (g\,\func{$\circ$}\,\func{F}_{n,m}))}{(\push{(x\,\con{,}\,y)})}$ unfolds to
\begin{align}
\label{eq:+join}
\congpt{f}{(\ap{\func{F}_{n,m}}{(\func{$\ell$}(x\,\con{,}\,y))})}\,\func{$\cdot$}\, \congpt{g}{(\ap{\func{F}_{n,m}}{(\func{$\ell$}(x\,\con{,}\,y))})}
\end{align}
Hence, comparing \eqref{eq:+susp} and \eqref{eq:+join}, we see that it is enough to show that $$\ap{\func{F}_{n,m}}{(\func{$\ell$}(x\,\con{,}\,y))} \,\func{$\equiv$}\, \func{$\sigma$}(x \,\func{$\smile$}\, y)$$
Unfolding $\func{$\ell$}$, we get
\begin{align*}
\ap{\func{F}_{n,m}}{(\func{$\ell$}(x\,\con{,}\,y))} \,&\func{$\equiv$}\,
\func{$\sigma$}(\cupprod{\star_{\sphere{n}}}{\star_{\sphere{m}}}) \,\func{$\cdot$}\,
\func{$\sigma$}(\cupprod{x}{\star_{\sphere{m}}})^{\func{$-1$}} \,\func{$\cdot$}\,
\func{$\sigma$}(\cupprod{x}{y})  \,\func{$\cdot$}\,
\func{$\sigma$}(\cupprod{\star_{\sphere{n}}}{y})^{\func{$-1$}}
\\
&\func{$\equiv$}\,
\func{$\sigma$}\,\north \,\func{$\cdot$}\,
\func{$\sigma$}\,\north^{\func{$-1$}} \,\func{$\cdot$}\,
\func{$\sigma$}(\cupprod{x}{y})  \,\func{$\cdot$}\,
\func{$\sigma$}\,\north^{\func{$-1$}}
\\ &\func{$\equiv$}\,\func{$\sigma$}(\cupprod{x}{y}) \qedhere
\end{align*}
\end{proof}

\begin{prop}
For any pointed type $A$, the set $\func{$\pi$}_{n+m+1}^*(A)$ is a
group with group structure induced by $\func{$+^*$}$. Furthermore,
pre-composition \textnormal{$(\func{F}_{n,m})^* :
  \func{$\pi$}_{n+m+1}(A) \to \func{$\pi$}_{n+m+1}^*(A)$} is an
isomorphism.
\end{prop}
\begin{proof}
We know that $(\func{F}_{n,m})^*$ is an equivalence of types. By
\autoref{prop:+-pres} and the Structure Identity
Principle~\cite[Section 9.8]{HoTT13}, it induces
a path
\[
(\func{$\pi$}_{n+m+1}(A),\func{$+^{\func{\textnormal{Susp}}}$})
\, \func{$\equiv$}\,(\func{$\pi$}^*_{n+m+1}(A), \func{$+^*$})
\]
of raw monoids (i.e.\ elements of type $\Sigma_{A : \Type}(A\,
\,\func{$\times$}\, A \to A)$). Since the the left-hand side of this
equality can be extended to form a group, so can the right-hand
side. This is precisely what we set out to show.
\end{proof}
The following result follows in exactly the same manner.
\begin{prop}
$\func{$\pi$}_{n+m+1}^*$ is functorial with its action on maps being
  defined by post-composition.
\end{prop}

\subsection{The new synthetic proof that $\tyIso{\func{$\pi$}_4(\sphere{3})}{\bZ}/2\bZ$}
\label{theproof}

Let us now return to the new proof. We will use $\cupprodop$ from
above in dimensions $\sphere{1}\,\func{$\times$}\,\sphere{1} \to
\sphere{2}$. We remark that, by~\autoref{prop:cup-comm}, it is
anti-commutative in these dimensions. In order to make the following
constructions somewhat more direct, let us return to the $\base/\Loop$
definition of $\sphere{1}$. Under the equivalence,
$\Susp{\func{Bool}}\,\func{$\simeq$}\,\sphere{1}$, the multiplication
is described by
\ExecuteMetaData[agda/latex/background.tex]{circprod}
In addition to anti-commutativity and associativity, we have the following distributivity-like fact about $\cupprodop$:
\begin{lem}\label{prop:circprod-props}
  For $x,y : \sphere{1}$, we have \textnormal{$\cupprodalt{x}{(x\,\func{+}\,y)}\,\func{$\equiv$}\,{\cupprodalt{x}{y}}$}
\end{lem}
\begin{proof}
  We proceed by $\sphere{1}$-induction on $x$. The equality $\cupprodalt{\base}{(\base\,\func{+}\,y)}\,\func{$\equiv$}\,{\cupprodalt{\base}{y}}$ holds by $\refl$, so we are left to verify the equality
  \[
  \ap{(x\mapsto \cupprod{x}{(x\,\func{+}\,y)})}{\Loop}\,\func{$\equiv$}\,\func{$\sigma$}\,y
  \]
  Simplifying the left-hand side using functoriality of binary $\func{cong}$~\cite[Definition 1]{LM24}, we get
  \begin{align*}
    \ap{(x\mapsto \cupprod{x}{(x\,\func{+}\,y)})}{\Loop}\,&\func{$\equiv$}\,\ap{(x\mapsto \cupprod{\north}{(x\,\func{+}\,y)})}{\Loop} \\
    \!\!&\phantom{\func{$\equiv$}\!\!}\func{$\cdot$}\, \ap{(x\mapsto \cupprod{x}{(\north\,\func{+}\,y)})}{\Loop}\\
    &:=\,\refl \,\func{$\cdot$} \, \func{$\sigma$}\,y \, \func{$\equiv$}\, \func{$\sigma$}\,y \qedhere
  \end{align*}
\end{proof}
We now redefine
$\tyIso{\func{$\pi$}_3{(\sphere{2})}}{\bZ}$ via the following
decomposition, primarily defined in terms of post- and pre-composition
with $\func{F}_{1,1} : \JoinS \,\func{$\cong$}\, \sphere{3} $ and its inverse. In what follows, let us simply write $\func{F} := \func{F}_{1,1}$ and $\func{$\pi$}_3^*(A) := \func{$\pi$}_{1+1+1}^{*}(A) := \truncT{0}{\JoinS \to_\star A}$. We also remind the reader of the map $\func{h} : \JoinS \to \sphere{2}$ from \autoref{def:hopfmap} for which $\func{h}_*$ is an isomorphism--this follows from~\autoref{prop:hopfcomp}.
\begin{defi}\label{def:newiso}
  Let $\theta : \tyIso{\func{$\pi$}_3{(\sphere{2})}}{\bZ}$ be defined
  by the following sequence of isomorphisms
  \[
\begin{tikzcd}[column sep=3em]
  \func{$\pi$}_3(\sphere{2}) \arrow["\textnormal{\func{F}}^*",r] &
  \func{$\pi$}^*_3(\sphere{2}) \arrow["(\textnormal{\func{h}}_*)^{-1}",r] &
  \func{$\pi$}^*_3(\JoinS) \arrow["\textnormal{\func{F}}_*",r] &
  \func{$\pi$}^*_3(\sphere{3}) \arrow["(\textnormal{\func{F}}^{-1})^{*}",r] &
  \func{$\pi$}_3(\sphere{3}) \arrow["\xi",r] &
  \bZ
\end{tikzcd}
\]
where the last map can be chosen to be any reasonable description of the
isomorphism $\xi : \tyIso{\func{$\pi$}_3(\sphere{3}) }{\bZ}$ sending
$i_3$ to $1$.
\end{defi}

The goal is to trace the image of
$[i_2,i_2] : \func{$\pi$}_3(\sphere{2})$ under $\theta$. Let us define
the following three underlying functions of elements
$\func{$\eta$}_1 : \func{$\pi$}^*_3(\sphere{2})$,
$\func{$\eta$}_2:\func{$\pi$}^*_3(\JoinS)$ and
$\func{$\eta$}_3:\func{$\pi$}^*_3(\sphere{3})$:
\ExecuteMetaData[agda/latex/background.tex]{eta1}
\ExecuteMetaData[agda/latex/background.tex]{eta2}
\ExecuteMetaData[agda/latex/background.tex]{eta3}

The claim is now that the image of $[i_2,i_2]$ under the chain of
isomorphisms can be described as follows:
\[
\begin{tikzcd}[column sep=4em]
  [i_2,i_2]\arrow["\func{F}^*",r,mapsto] &
  \func{$\eta$}_1 \arrow["(\func{h}_*)^{-1}",r,mapsto] &
  \func{$\eta$}_2 \arrow["\func{F}_*",r] &
  \func{$\eta$}_3 \arrow["(\func{F}^{-1})^{*}",r,mapsto] &
  (-2)i_3 \arrow["\xi",r,mapsto] &
  \pm 2
\end{tikzcd}
\]

\begin{lem}\label{lem:trace1}
  \textnormal{$\tyPath{\func{F}^*\,[i_2,i_2]}{\func{$\eta$}_1}$}
\end{lem}
\begin{proof}
  The definition of $\func{$\eta$}_1$ matches that of $\trunc{\Nabla
    \circ \func{W}} : \func{$\pi$}_3^{*}(\sphere{2})$, and so the
  statement holds by construction of the Whitehead product.
\end{proof}

\begin{lem}\label{lem:trace2}
  \textnormal{$\tyPath{(\func{h}_*)^{-1}\,\func{$\eta$}_1}{\func{$\eta$}_2}$}
\end{lem}
\begin{proof}
  Applying $\func{h}_*$ on both sides gives the equation
  \textnormal{$\tyPath{\func{$\eta$}_1}{\func{h}_*\,\func{$\eta$}_2}$}. Thus,
  we are done if we can show that $\etafun{1}\,a
  \,\func{$\equiv$}\,\func{h}\,(\etafun{2}\,a)$ for $a : \JoinS$. We
  do it by induction on $a$. When $a$ is $\inl{x}$ or $\inr{y}$, the
  equality holds by $\refl$. Thus, it remains to show that
  \[
  \ap{\etafun{1}}{(\push{(x\,\con{,}\,y)})} \,\func{$\equiv$}\,
  \ap{(\func{h}\,\func{$\circ$}\,\etafun{2})}{(\push{(x\,\con{,}\,y)})}
  \]
  We show the identity by unfolding the right-hand side:
  \begin{align*}
    \ap{(\func{h}\,\func{$\circ$}\,\etafun{2})}{(\push{(x\,\con{,}\,y)})}\,&{:=}\, \ap{\func{h}}{(\push{(y\,\func{$-$}\,x \,\con{,}\, \func{$-$}\,x)}^{\func{$-1$}} \,\func{$\cdot$}\, \push{(y\,\func{$-$}\,x \,\con{,}\, y)})}  \\
    &\func{$\equiv$}\, \ap{\func{h}}{(\push{(y\,\func{$-$}\,x \,\con{,}\, \func{$-$}\,x)})^{\func{$-1$}}} \,\func{$\cdot$}\, \ap{\func{h}}{(\push{(y\,\func{$-$}\,x \,\con{,}\, y)})}\\
    &{:=}\, \func{$\sigma$}((\func{$-$}\,x) \,\func{$-$}\,(y\,\func{$-$}\,x))^{\func{$-1$}}\,\func{$\cdot$}\,\func{$\sigma$}(y \,\func{$-$}\,(y\,\func{$-$}\,x)) \\
    &\func{$\equiv$}\, \func{$\sigma$}(\func{$-$}\,y)^{\func{$-1$}}\,\func{$\cdot$}\,\func{$\sigma$}\,x \\
    &\func{$\equiv$} \,\func{$\sigma$}\,y\,\func{$\cdot$}\,\func{$\sigma$}\,x\\
    &{=:}\, \ap{\etafun{1}}{(\push{(x\,\con{,}\,y)})} \qedhere
  \end{align*}
\end{proof}

\begin{lem}\label{lem:trace3}
  \textnormal{$\tyPath{\func{F}_{*}\,\func{$\eta$}_2}{\func{$\eta$}_3}$}
\end{lem}
\begin{proof}
  The identity follows if we can show that $\func{F} (\etafun{2}\,a) \,\func{$\equiv$} \,\etafun{3}\,a$ for $a : \JoinS$. Again, the identity holds by $\refl$ when $a$ is $\inl{x}$ or $\inr{y}$. So it remains to show that
  \[\ap{(\func{F} \,\circ\,\etafun{2})}{(\push{(x\,\con{,}\,y)})} \,\func{$\equiv$}\, \ap{\etafun{3}}{(\push{(x\,\con{,}\,y)})} \]
  Just like in the proof of~\autoref{lem:trace3}, we show this simply by unfolding the definitions of the, in this case, left-hand side. We get:
  \begin{align*}
    \ap{(\func{F}\,\func{$\circ$}\,\etafun{2})}{(\push{(x\,\con{,}\,y)})}\,&{:=}\, \ap{\func{F}}{(\push{(y\,\func{$-$}\,x \,\con{,}\, \func{$-$}\,x)}^{\func{$-1$}} \,\func{$\cdot$}\, \push{(y\,\func{$-$}\,x \,\con{,}\, y)})}  \\
    &\func{$\equiv$}\, \ap{\func{F}}{(\push{(y\,\func{$-$}\,x \,\con{,}\, \func{$-$}\,x)})^{\func{$-1$}}} \,\func{$\cdot$}\, \ap{\func{F}}{(\push{(y\,\func{$-$}\,x \,\con{,}\, y)})}\\
    &{:=}\, \func{$\sigma$}(\cupprod{(y\,\func{$-$}\,x)}{(\func{$-$}\,x)})^{\func{$-1$}}\,\func{$\cdot$}\,\func{$\sigma$}(\cupprod{(y\,\func{$-$}\,x)}{y}) \\
    &\func{$\equiv$}\, \func{$\sigma$}(\cupprod{(\func{$-$}\,x)}{((\func{$-$}\,x)\,\func{$+$}\,y)})\,\func{$\cdot$}\,\func{$\sigma$}(\cupprod{y}{(y\,\func{$-$}\,x)})^{\func{$-1$}} \\
    &\func{$\equiv$} \,\func{$\sigma$}(\cupprod{(\func{$-$}\,x)}{y})\,\func{$\cdot$}\,\func{$\sigma$}(\cupprod{y}{(\func{$-$}\,x)})^{\func{$-1$}}\\
    &\func{$\equiv$} \,\func{$\sigma$}(\cupprod{x}{y})^{\func{$-1$}}\,\func{$\cdot$}\,\func{$\sigma$}(\cupprod{y}{(\func{$-$}\,x)})^{\func{$-1$}}\\
    &\func{$\equiv$} \,\func{$\sigma$}(\cupprod{x}{y})^{\func{$-1$}}\,\func{$\cdot$}\,\func{$\sigma$}(\cupprod{x}{y})^{\func{$-1$}}\\
    &{=:}\, \ap{\etafun{3}}{(\push{(x\,\con{,}\,y)})}
  \end{align*}
where the fourth and seventh equalities come from anti-commutativity and the fifth
equality from \autoref{prop:circprod-props}. The fact that $\func{$\sigma$}$ commutes with inversion is used throughout.
\end{proof}

\begin{thm}
  $\tyIso{\func{$\pi$}_4(\sphere{3})}{\bZ}/2\bZ$
\end{thm}
\begin{proof}
  By uniqueness (up to a sign) of isomorphisms
  $\tyIso{\func{$\pi$}_3(\sphere{2})}{\bZ}$, it suffices, according to
  \autoref{cor:main}, to show that the image of $[i_2,i_2]$ under
  $\theta$ is $\pm 2$. That is:
  \[\tyPath{(\xi \circ (\func{F}^{-1})^* \circ \func{F}_* \circ (\func{h}_*)^{-1} \circ \func{F}^*) [i_2,i_2]}{\pm2} \]
  By~\autoref{lem:trace1}, \autoref{lem:trace2} and
  \autoref{lem:trace3}, it suffices to show that
  \[\tyPath{(\xi \circ (\func{F}^{-1})^*)\,\func{$\eta$}_3}{\pm 2}\]
  One can easily show that
  $\tyPath{\func{F}^{-1}\,\func{$\eta$}_3}{(-2)\,i_3}$, and hence
  \[\tyPath{(\xi \circ (\func{F}^{-1})^*)\,\func{$\eta$}_3}{(-2)\,(\xi\,i_3)}\,\func{$\equiv$}\,-2 \qedhere \]
\end{proof}

In addition to providing a much shorter proof of
$\tyIso{\func{$\pi$}_4(\sphere{3})}{\bZ}/2\bZ$, this gives us a
sequence of new Brunerie numbers,
$\brunerie_1,\brunerie_2,\brunerie_3 : \bZ$, of decreasing complexity:
\begin{align*}
  \brunerie_1 &= (\xi \circ (\func{F}^{-1})^* \circ \func{F}_* \circ (\func{h}_*)^{-1})\,\func{$\eta$}_1 \\
  \brunerie_2 &= (\xi \circ (\func{F}^{-1})^* \circ \func{F}_*)\,\func{$\eta$}_2 \\
  \brunerie_3 &= (\xi \circ (\func{F}^{-1})^*)\,\func{$\eta$}_3
\end{align*}
This gives new hope for Brunerie's conjecture about a proof by
normalisation. This may be captured as follows:
\begin{thm}[New Brunerie numbers]
  If either of $\brunerie_1,\brunerie_2,\brunerie_3 : \bZ$ normalises
  to $\pm 2$, then $\tyIso{\func{$\pi$}_4(\sphere{3})}{\bZ}/2\bZ$.
\end{thm}
Ideally, we could normalise $\brunerie_1$. This, however, turns out to
be difficult, as it does not bypass the main hurdle of computing the
inverse of the isomorphism
$\tyIso{\func{$\pi$}^*_3(\sphere{2})}{\func{$\pi$}_3^*(\JoinS)}$
induced by the Hopf map, which has a rather indirect construction
coming from the LES of homotopy groups associated to the Hopf
fibration. This problem does not apply to $\brunerie_2$, for which the
computation does not rely on the problematic inverse. Unfortunately,
also $\brunerie_2$ fails to normalise in reasonable time in
\CubicalAgda.  This is surprising, as the only maps playing a
fundamental role here are two applications of the equivalence
$\tyEquiv{\JoinS}{\sphere{3}}$, which is not too involved, and one
application of $\xi$ which may be compactly described via
\[\func{$\pi$}_3{(\sphere{3})}\xrightarrow{\trunc{\!\con{\_}\!}_*} \cohom{3}{\sphere{3}}  \xrightarrow{\func{$\cong$}} \bZ \]
and computes relatively well if the last isomorphism is constructed as
in~\cite{BLM22}.\footnote{As noted in~\cite{LjungstromMsc20}, the
  Freudenthal suspension theorem should be avoided here as it has a
  tendency to lead to very slow computations. This is another way in
  which we deviate from Brunerie's $\brunerie$.} We have hence, at
the time of writing, not been able to normalise even $\brunerie_2$,
despite many optimisations of the functions involved. We are, however,
able to normalise $\brunerie_3$ after some minor modifications to
$\func{$\eta$}_3$ and the map $\func{$\pi$}_3^*(\sphere{3}) \to
\bZ$. This optimised version of $\brunerie_3$, normalises to $-2$ in
\CubicalAgda in just under $4$ seconds, thereby giving us an at least
partially computer-assisted proof of
$\tyIso{\func{$\pi$}_4(\sphere{3})}{\bZ}/2\bZ$.

We emphasise again that $\brunerie_2$ is a vastly simplified version
of $\func{$\beta$}$ since the isomorphism
$\tyIso{\func{$\pi$}_3(\sphere{2})}{\func{$\pi$}_3(\sphere{3})}$ never
has to be computed. Hence, it is rather surprising that computations
break down already at this stage. This tells us that \CubicalAgda has
a long way to go before any direct computation of the original
$\brunerie$ is feasible. We hope that this could be useful for
benchmarking in future optimisations of \CubicalAgda and related
systems.

Finally, we address the elephant in the room: why is there a minus
sign popping up? In other words, have we really chosen the, in some
way, canonical isomorphism? The isomorphism
$\tyIso{\func{$\pi$}_3(\sphere{3})}{\bZ}$ maps, as expected, $i_3$ to
$1$, so it can hardly be the culprit. Neither can the equivalence
$\textnormal{\func{F}}:\tyEquiv{\JoinS}{\sphere{3}}$, since it is
applied equally in the constructions of $\hopf$ and of $[i_2,i_2]$. We
could, however, have defined the $\push{}$\!-case for $\func{h}$ by
\ExecuteMetaData[agda/latex/background.tex]{alt-hopf}
in which case $\theta$ would have sent $[i_2,i_2]$ to $2$ and
$\hopf$ to $1$ (note that this is only possible since altering
$\func{h}$ would alter the definition of $\theta$). The construction
of $\func{h}$ that we have given is, however, precisely the one which
fell out by unfolding our formalisation Brunerie's construction of the
corresponding map. If this indeed is what Brunerie intended, we may
also conclude that the original Brunerie number $\brunerie$ is equal
to $-2$. We stress that this merely is a fun fact and of no
mathematical importance to Brunerie's proof or our formalisation.

\subsection{A stand-alone proof of Brunerie's theorem?}
\label{standalone}

We saw above that the new proof of $\func{$\beta$}
\,\func{$\equiv$}\, \pm 2$ together with \autoref{cor:main} implies
Brunerie's theorem. However, what conclusions can we draw concerning
the cardinality of $\func{$\pi$}_4(\sphere{3})$ in the absence of
\autoref{cor:main}? In other words, how self-contained is the new
proof? While the fact that $\func{$\beta$} \,\func{$\equiv$}\, \pm 2$
does not automatically imply that
$\func{$\pi$}_4(\sphere{3})\,\func{$\cong$}\,\bZ/2\bZ$, it does
provide all ingredients necessary for a stand-alone proof of the
following fact:
\begin{thm}\label{thm:self-contained}
  If $\func{$\pi$}_4(\sphere{3}) \,\func{$\not\simeq$}\, \func{$\mathbbm{1}$}$, then $\func{$\pi$}_4(\sphere{3})\,\func{$\cong$}\,\bZ/2\bZ$.
\end{thm}
Before we prove \autoref{thm:self-contained}, we need to analyse the
action of suspension on Whitehead products and, in particular, on $[i_2,i_2] :
\func{$\pi$}_3(\sphere{2})$. In what follows, let $A,B$ and $C$ be pointed types and let us fix two pointed functions $f:\Susp{A} \to_\star C$ and $g : \Susp{B} \to_\star C$. The Whitehead product of $f$ and $g$ can be understood as the composition
\[\JJoin{A\,}{B} \xrightarrow{\, \func{W} \,} \Susp{A}\,\func{$\vee$}\,\Susp{B} \xrightarrow{f\,\func{$\vee$}\,g} C\,\func{$\vee$}\,C \xrightarrow{\,\func{$\nabla$}\,} C\]
We remark that this construction has been independently studied by Cagne et al.~\cite[Definition 6.3]{cagne2024symmetries} who call it the `generalised Whitehead product'.
After a bit of massaging, this function can be given a very simple description:
\ExecuteMetaData[agda/latex/background.tex]{wh}
We remark that this composition gives $\etafun{1}$ when $A = B = \sphere{1}$, $C = \sphere{2}$ and $f = g = \idfun_{\sphere{2}}$.

Our aim is to show that $\fdotg$ vanishes under suspension. To this end, let us consider a function very similar to $\fdotg$:
\ExecuteMetaData[agda/latex/background.tex]{gamma}
Despite the similarity of $\fdotg$ and $\func{$\gamma$}$, the latter turns out to be trivial.
\begin{lem}\label{lem:gamma-triv}
  \func{$\gamma$} is constant.
\end{lem}
\begin{proof}
  We show that $\func{$\gamma$}\,a \,\func{$\equiv$}\,\star_C$ for all
  $a : \JoinS$ by induction on $a$. When $a$ is $\inl{x}$, the
  left-hand side reduces to $\star_C$, so we need to provide a path
  $\star_C\,\func{$\equiv$}\,\star_C$. Instead of choosing the obvious
  path $\refl$, we provide $\congpt{f}{(\func{$\sigma$}\, x)} :
  \star_C\,\func{$\equiv$}\,\star_C$. When $a$ is $\inr{y}$, we have
  the same goal. This time, we provide the path
  $\congpt{g}{(\func{$\sigma$}\, y)}^{\func{$-1$}}$. For the final
  step, i.e.\ the action of $\func{$\gamma$}$ on
  $\push{(x\,\con{,}\,y)}$, we need to provide a filler of the
  following square of paths:
  \[
\begin{tikzcd}[ampersand replacement=\&,column sep=5em]
	\star_C \&\& \star_C \\
	\star_C \&\& \star_C
	\arrow["\refl", Rightarrow, no head, from=1-1, to=1-3]
	\arrow["{\congpt{f}{(\func{$\sigma$}\, x)}}", from=2-1, to=1-1]
	\arrow["{\congpt{f}{(\func{$\sigma$}\, x)}\,\func{$\cdot$}\,\congpt{g}{(\func{$\sigma$}\, y)}}"', from=2-1, to=2-3]
	\arrow["{(\congpt{g}{(\func{$\sigma$}\, y)})^{\func{$-1$}}}"', from=2-3, to=1-3]
\end{tikzcd}
\]
Squares of this shape always have a filler by definition of path
composition, and thus the statement holds.
\end{proof}
Now, although $\fdotg$ and $\func{$\gamma$}$ may look similar, it it
is now, in light of~\autoref{lem:gamma-triv}, clear that they are not
the same. This happens because the actions of the functions on
$\push{(x\,\con{,}\,y)}$ only are the same up to commutation of
paths---something which is not always legal in the possibly
non-commutative loop space $\LoopSpace{C}$. Nevertheless, after
suspending the function, the situation is different:
\begin{lem}\label{lem:wh-vanish}
  The pointed functions $\Susp{(\fdotg)}, \Susp{\func{$\gamma$}} :
  \Susp{(\JJoin{A\,}{B})} \to_\star \Susp{C}$ are equal.
\end{lem}
\begin{proof}
  Under the adjunction $\func{Susp} \dashv \func{$\Omega$}$, it is
  enough to show that for every $a : \JJoin{A\,}{B}$, we have an
  equality of loops in ${\LoopSpace{(\Susp{C})}}$:
  \[
  \func{$\sigma$}((\fdotg)\,a) \,\func{$\equiv$} \func{$\sigma$}(\func{$\gamma$}\,a) 
  \]
  We proceed by induction on $a$. When $a$ is $\inl{x}$ or $\inr{y}$,
  the equality holds by $\refl$. Thus, it remains to show that
  \[
  \ap{\func{$\sigma$}}{(\ap{(\fdotg)}{(\push{(x\,\con{,}\,y)})})} \,\func{$\equiv$}\, \ap{\func{$\sigma$}}{(\ap{\func{$\gamma$}}{(\push{(x\,\con{,}\,y)})})}
  \]
  As before, this is a simple exercise in unfolding the definitions of
  each respective function:
  \begin{align*}
    \ap{\func{$\sigma$}}{(\ap{(\fdotg)}{(\push{(x\,\con{,}\,y)})})} \,&{:=}\, \ap{\func{$\sigma$}}{(\congpt{g}{(\func{$\sigma$}\,y)}\,\func{$\cdot$}\, \congpt{f}{(\func{$\sigma$}\,x)})}\\
    &\func{$\equiv$}\, \ap{\func{$\sigma$}}{(\congpt{g}{(\func{$\sigma$}\,y)})}\,\func{$\cdot$}\, \ap{\func{$\sigma$}}{(\congpt{f}{(\func{$\sigma$}\,x)})}\\
    &\func{$\equiv$}\, \ap{\func{$\sigma$}}{(\congpt{f}{(\func{$\sigma$}\,x)})}\,\func{$\cdot$}\, \ap{\func{$\sigma$}}{(\congpt{g}{(\func{$\sigma$}\,y)})}  \tag{EH}\label{eq:eh-appl} \\
    &\func{$\equiv$}\,\ap{\func{$\sigma$}}{(\congpt{f}{(\func{$\sigma$}\,x)}\,\func{$\cdot$}\, \congpt{g}{(\func{$\sigma$}\,y)})}\\
    &\func{$\equiv$}\,\ap{\func{$\gamma$}}{(\push{(x\,\con{,}\,y)})}
  \end{align*}
where the step labelled \eqref{eq:eh-appl} is an application of the
Eckmann-Hilton argument which says that path composition in
$\LoopSpaceN{2}{A}$ is commutative for any pointed type
$A$~\cite[Theorem 2.1.6]{HoTT13}. In particular, since we may
interpret $\ap{\func{$\sigma$}}{(\congpt{f}{(\func{$\sigma$}\,x)})}$
and $\ap{\func{$\sigma$}}{(\congpt{g}{(\func{$\sigma$}\,y)})}$ as
loops in $\LoopSpaceN{2}{(\Susp{C})}$, the identity holds.
\end{proof}
\begin{prop}\label{prop:wh-vanish}
  The pointed function $\Susp{(\fdotg)} : \Susp{(\JJoin{A\,}{B})} \to_\star \Susp{C}$ is constant.
\end{prop}
\begin{proof}
  Since $\func{$\gamma$}$ is constant and constant functions are preserved by suspension, \autoref{lem:gamma-triv} gives us the desired equality of (pointed) functions:
  \[\Susp{(\fdotg)} \,\func{$\equiv$} \, \Susp{\func{$\gamma$}} \,\func{$\equiv$}\, \func{const}_{C} \qedhere \]
\end{proof}
As we have seen before, setting $A = \sphere{n}$ and $B=\sphere{m}$ in
the definition of $\fdotg$, so that $f : \sphere{n+1} \to_\star C$ and
$g : \sphere{m+1} \to_\star C$ we obtain the usual Whitehead
product. $[\trunc{f},\trunc{g}] : \func{$\pi$}_{n+m+1}(C)$, that is
\[
[\trunc{f},\trunc{g}] \,\func{$\equiv$}\, \trunc{(\fdotg)\,\func{$\circ$}\,\func{F}_{n,m}^{-1}}
\]
Let us translate \autoref{prop:wh-vanish} to a result concerning these maps.
\begin{prop}
  For $f : \sphere{n} \to_\star C$ and $g : \sphere{m} \to_\star C$,
  their Whitehead product
  \textnormal{$(\fdotg)\,\func{$\circ$}\,\func{{F}}_{n,m}^{-1}$}
  vanishes under suspension,
  i.e.\ \textnormal{\[\Susp{(\trunc{(\fdotg)\,\func{$\circ$}\,\func{F}_{n,m}^{-1}})}
    \,\func{$\equiv$}\, \func{{const}}_{C} \]}
\end{prop}
\begin{proof}
  The result follows immediately from the fact that the action of
  suspension $\func{Susp}:(X \to_\star Y) \to (\Susp{X} \to_\star
  \Susp{Y})$ is functorial and from \autoref{prop:wh-vanish}.
\end{proof}
We get the following classically well-known theorem as an immediate
corollary:
\begin{thm}\label{tmm:wh-vanish}
  For any $x: \func{$\pi$}_{n+1}(C)$ and $y : \func{$\pi$}_{m+1}(C)$,
  the Whitehead product $[x,y] : \func{$\pi$}_{n+m+1}(C)$ lies in the
  kernel of the suspension map $\func{$\sigma$}_* :
  \func{$\pi$}_{n+m+1}(C) \to \func{$\pi$}_{n+m+2}(\Susp{C})$
\end{thm}
We now have all that we need in order to
prove~\autoref{thm:self-contained}.
\begin{proof}[Proof of \autoref{thm:self-contained}]
  By the Freudenthal suspension theorem, we know that
  $\func{$\sigma$}_* : \func{$\pi$}_3(\sphere{2}) \to
  \func{$\pi$}_4(\sphere{3})$ is surjective. Furthermore, we know that
  the domain of this function is isomorphic to $\bZ$ via $\theta$ from
  \autoref{def:newiso} and thus we have a surjection
  $\func{$\sigma$}_* \,\func{$\circ$}\,\theta^{-1} : \bZ \to
  \func{$\pi$}_4(\sphere{3})$. We know from the new direct calculation
  of the Brunerie number that $\theta^{-1}(-2) \,\func{$\equiv$}\,
  [i_2,i_2]$ and thus we have
  \[
  \func{$\sigma$}_* (\theta^{-1}(-2)) \,\func{$\equiv$}\,\func{$\sigma$}_*[i_2,i_2] \,\func{$\equiv$}\, 0_{\func{$\pi$}_4(\sphere{3})}
  \]
where the second equality comes from \autoref{tmm:wh-vanish}. Hence,
we have shown that there exists a surjection from $\bZ$ onto
$\func{$\pi$}_4(\sphere{3})$ with $-2$ in its kernel. This implies the
theorem.
\end{proof}
Now, with \autoref{thm:self-contained} in mind, we seem to be very
close to having produced a remarkably short proof of Brunerie's
theorem. All that remains is showing that $\func{$\pi$}_4(\sphere{3})$
is not trivial. This, however, turns out not to be entirely
straightforward. One possible proof uses the so called \emph{Steenrod
  Squares}. This is a cohomology operation which was originally
defined in HoTT by Brunerie~\cite{Brunerie17} and whose theory was
recently made available in HoTT by
Ljungström~and~Wärn~\cite{SteenrodSquares}. Such an approach, however,
can hardly be said to simplify Brunerie's original proof, as the
Steenrod Squares are rather advanced constructions.
A solution to this problem which would truly be impressive would be a
direct construction of two elements $x,y : \func{$\pi$}_4(\sphere{3})$
and a proof that $x \,\func{$\not\equiv$}\,y$. While this appears to
be difficult to do by hand, we can, since we are working
constructively, reformulate this problem as a computational challenge.
\begin{challenge}
  Construct a function \textnormal{$f : \func{$\pi$}_4(\sphere{3})\to
    \func{{Bool}}$} and an element $e : \func{$\pi$}_4(\sphere{3})$
  such that
  \begin{itemize}
  \item $f\,0_{\func{$\pi$}_4(\sphere{3})}$ computes to
    \textnormal{$\truebool$} and
    \item $f\,e$ computes to \textnormal{$\falsebool$}.
  \end{itemize}
\end{challenge}
In fact, such a computation was successfully run by Jack~\cite{Jack}
in \cubicaltt \cite{ctt}. Unfortunately, \CubicalAgda has not yet been
able to perform the computation.

%% file: sections/conclusion.tex
\section{Conclusion}
\label{sec:conclusion}

In this paper, we have presented three formalisations of
$\tyIso{\func{$\pi$}_4(\sphere{3})}{\bZ}/2\bZ$ in the \CubicalAgda
system. For the different proofs that
$\abs{\brunerie}\,\func{$\equiv$}\,2$, the line count is roughly as
follows:
\begin{enumerate}
\item Brunerie's original proof [$\sim9,000$ LOC]
\item A direct calculation of $\brunerie$ [$\sim600$ LOC]
\item A computer-assisted reformulation of (2) [$\sim400$ LOC]
\end{enumerate}
As always, the number of lines of code (LOC) should be taken with a
grain of salt. First, the $9,000$ LOC in the first formalisation
exclude over $8,000$ LOC
from~\cite{BMFormalization,FreudenthalFormalization,BLM22} which we
have imported as libraries. In addition, these numbers also exclude
many elementary results used in the formalisation, including $\sim
9000$ LOC for Chapters 1--3. We also stress that the line count for
formalisations (2) and (3) only concern the part of the proof discussed
in~\autoref{sec:computation}.

Formalisation (1), which constituted the bulk of this paper, was a
formalisation of Brunerie's pen-and-paper proof, taking some
convenient shortcuts when possible. The problem of formalising
Brunerie's proof has been a widely discussed open problem in HoTT/UF,
and we hope that our efforts here provide a satisfactory solution to
it. Formalisations (2) and (3) were of a simplified calculation
of the Brunerie number, $\brunerie$. The very similar proofs (2) and (3) differ in
that (3) uses \CubicalAgda to carry out part of the computation of the
new Brunerie number automatically. Perhaps equally important, we have
seen that (3) provides us with new Brunerie numbers
$\brunerie_1, \brunerie_2 : \bZ$ which are far simpler than the
original one, but still do not normalise in a reasonable amount of
time. Our hope is that these can prove useful in future optimisations
of \CubicalAgda and related systems, as they could help shed some
light on where the normalisation of the original Brunerie number
breaks down.

We remark that proofs (1) and (2) could be done in Book HoTT and do
not use any cubical machinery in a fundamental way, making them
interpretable in any suitably structured $(\infty,1)$-topos
\cite{Shulman19}. We hence claim that, in our formalisations, we do
not crucially rely on computations using univalence and HITs to prove
anything that we could not have proved by hand in Book
HoTT. Nevertheless, the \CubicalAgda system has been very helpful in
the formalisation, primarily due to its native support for HITs and
definitional computation rules for higher constructors. Formalisation
(3), however, is only valid in a system with computational support for
univalence as it crucially relies on normalisation of proof terms
involving univalence. It would be interesting to run this in other
cubical systems, like \cubicaltt \cite{ctt}, \redtt{} \cite{redtt},
\cooltt{} \cite{cooltt}, etc.

In addition to the above, we have also taken the opportunity to
include some important constructions and results concerning joins of
spheres and Whitehead products. In particular, we have given a very
explicit definition of the decomposition of spheres into joins of
spheres, given a new construction of homotopy groups in terms of maps
out of joins of spheres and shown that Whitehead products vanish under
suspension. The vanishing of Whitehead products allowed us to extend
(2) to a stand-alone proof of the fact that
$\func{$\pi$}_4(\sphere{3})$ is either trivial or isomorphic to
$\bZ/2\bZ$. Interestingly, another direct proof using an entirely
different approach of this very fact was recently announced by
Baker~\cite{BakerHoTTUF}. Baker's argument is concerned with showing
that a certain path constructed via the Eckmann-Hilton argument
generates $\func{$\pi$}_3(\sphere{2})$ and then concludes that two
times this generator must vanish under suspension due to the so called
\emph{syllepsis}~\cite{Syllepsis}. We leave it to future work to
investigate if anything interesting can be said about the relation
between Baker's proof and ours.

We also remark that our formalisation of Brunerie's proof does not
cover all results of Brunerie's thesis in full generality. For
instance, we have not developed his proof concerning Whitehead
products in full generality.  We leave this generalisation for future
work. This would tie in nicely with another possible direction of
future research, namely that of investigating whether the approach
outlined in \autoref{sec:computation} can be used to compute other
Whitehead products. In addition, describing their graded quasi-Lie
algebra structure is work in progress.

%% file: paper.bbl
\begin{thebibliography}{CCHM18}

\bibitem[{Agd}24]{Agda}
The {Agda Development Team}.
\newblock {The Agda Programming Language}, 2024.
\newblock URL: \url{http://wiki.portal.chalmers.se/agda/}.

\bibitem[AW09]{AwodeyWarren09}
Steve Awodey and Michael~A. Warren.
\newblock Homotopy theoretic models of identity types.
\newblock {\em Mathematical Proceedings of the Cambridge Philosophical
  Society}, 146(1):45--55, January 2009.
\newblock \href {https://doi.org/10.1017/S0305004108001783}
  {\path{doi:10.1017/S0305004108001783}}.

\bibitem[Bak24]{BakerHoTTUF}
Raymond Baker.
\newblock {Eckmann-Hilton and the Hopf Fibration}.
\newblock Extended abstract at \emph{Workshop on Homotopy Type Theory /
  Univalent Foundations (HoTT/UF24}, 2024.
\newblock URL:
  \url{https://hott-uf.github.io/2024/abstracts/HoTTUF_2024_paper_24.pdf}.

\bibitem[BHF18]{BuchholtzFavonia18}
Ulrik Buchholtz and Kuen-Bang Hou~Favonia.
\newblock {Cellular Cohomology in Homotopy Type Theory}.
\newblock In {\em Proceedings of the 33rd Annual ACM/IEEE Symposium on Logic in
  Computer Science}, LICS '18, pages 521--529, New York, NY, USA, 2018.
  Association for Computing Machinery.
\newblock \href {https://doi.org/10.1145/3209108.3209188}
  {\path{doi:10.1145/3209108.3209188}}.

\bibitem[BLM22]{BLM22}
Guillaume Brunerie, Axel Ljungstr\"{o}m, and Anders M\"{o}rtberg.
\newblock {Synthetic Integral Cohomology in Cubical Agda}.
\newblock In Florin Manea and Alex Simpson, editors, {\em 30th EACSL Annual
  Conference on Computer Science Logic (CSL 2022)}, volume 216 of {\em Leibniz
  International Proceedings in Informatics (LIPIcs)}, pages 11:1--11:19,
  Dagstuhl, Germany, 2022. Schloss Dagstuhl -- Leibniz-Zentrum f{\"u}r
  Informatik.
\newblock URL: \url{https://drops.dagstuhl.de/opus/volltexte/2022/15731}, \href
  {https://doi.org/10.4230/LIPIcs.CSL.2022.11}
  {\path{doi:10.4230/LIPIcs.CSL.2022.11}}.

\bibitem[Bru16a]{Brunerie16}
Guillaume Brunerie.
\newblock {\em On the homotopy groups of spheres in homotopy type theory}.
\newblock PhD thesis, Universit{\'e} Nice Sophia Antipolis, 2016.
\newblock URL: \url{http://arxiv.org/abs/1606.05916}.

\bibitem[Bru16b]{Brunerie17}
Guillaume Brunerie.
\newblock The steenrod squares in homotopy type theory.
\newblock Abstract at \emph{23rd International Conference on Types for Proofs
  and Programs ({TYPES} 2017)}, 2016.
\newblock URL: \url{https://types2017.elte.hu/proc.pdf#page=45}.

\bibitem[Bru18]{brunerie18}
Guillaume Brunerie.
\newblock Computer-generated proofs for the monoidal structure of the smash
  product.
\newblock \emph{Homotopy Type Theory Electronic Seminar Talks}, November 2018.
\newblock URL:
  \url{https://www.uwo.ca/math/faculty/kapulkin/seminars/hottest.html}.

\bibitem[Bru19]{GB-James}
Guillaume Brunerie.
\newblock {The James Construction and
  {\(\pi\)}\({}_{\mbox{4}}\)(S\({}^{\mbox{3}}\)) in Homotopy Type Theory}.
\newblock {\em Journal of Automated Reasoning}, 63:255--284, 2019.

\bibitem[Cav20]{FreudenthalFormalization}
Evan Cavallo.
\newblock {Formalisation of the Freudenthal Suspension Theorem}, 2020.
\newblock URL:
  \url{https://github.com/agda/cubical/blob/master/Cubical/Homotopy/Freudenthal.agda}.

\bibitem[CBKB24]{cagne2024symmetries}
Pierre Cagne, Ulrik Buchholtz, Nicolai Kraus, and Marc Bezem.
\newblock On symmetries of spheres in univalent foundations, 2024.
\newblock \href {https://arxiv.org/abs/2401.15037} {\path{arXiv:2401.15037}}.

\bibitem[CCHM]{ctt}
Cyril Cohen, Thierry Coquand, Simon Huber, and Anders M\"ortberg.
\newblock {\sc cubicaltt}: {Cubical Type Theory}.
\newblock Implementation available at
  \url{https://github.com/mortberg/cubicaltt}.

\bibitem[CCHM18]{CCHM18}
Cyril Cohen, Thierry Coquand, Simon Huber, and Anders M{\"o}rtberg.
\newblock {Cubical Type Theory: A Constructive Interpretation of the Univalence
  Axiom}.
\newblock In Tarmo Uustalu, editor, {\em 21st International Conference on Types
  for Proofs and Programs (TYPES 2015)}, volume~69 of {\em Leibniz
  International Proceedings in Informatics (LIPIcs)}, pages 5:1--5:34,
  Dagstuhl, Germany, 2018. Schloss Dagstuhl--Leibniz-Zentrum fuer Informatik.
\newblock \href {https://doi.org/10.4230/LIPIcs.TYPES.2015.5}
  {\path{doi:10.4230/LIPIcs.TYPES.2015.5}}.

\bibitem[CH19]{CavalloHarper19}
Evan Cavallo and Robert Harper.
\newblock {Higher Inductive Types in Cubical Computational Type Theory}.
\newblock {\em Proceedings of the ACM on Programming Languages},
  3(POPL):1:1--1:27, January 2019.
\newblock \href {https://doi.org/10.1145/3290314} {\path{doi:10.1145/3290314}}.

\bibitem[CHM18]{CoquandHuberMortberg18}
Thierry Coquand, Simon Huber, and Anders M\"{o}rtberg.
\newblock {On Higher Inductive Types in Cubical Type Theory}.
\newblock In {\em Proceedings of the 33rd Annual ACM/IEEE Symposium on Logic in
  Computer Science}, LICS '18, pages 255--264. ACM, 2018.
\newblock \href {https://doi.org/10.1145/3209108.3209197}
  {\path{doi:10.1145/3209108.3209197}}.

\bibitem[CS20]{christensen2020hurewicz}
J.~Daniel Christensen and Luis Scoccola.
\newblock {The Hurewicz theorem in Homotopy Type Theory}, 2020.
\newblock Preprint.
\newblock URL: \url{https://arxiv.org/abs/2007.05833}, \href
  {https://arxiv.org/abs/2007.05833} {\path{arXiv:2007.05833}}.

\bibitem[Hat02]{Hatcher2002}
Allen Hatcher.
\newblock {\em Algebraic Topology}.
\newblock Cambridge University Press, 2002.
\newblock URL: \url{https://pi.math.cornell.edu/~hatcher/AT/AT.pdf}.

\bibitem[HFLL16]{FavoniaFinster+16}
Kuen-Bang {Hou (Favonia)}, Eric Finster, Daniel~R. Licata, and Peter~LeFanu
  Lumsdaine.
\newblock {A Mechanization of the {Blakers-Massey} Connectivity Theorem in
  Homotopy Type Theory}.
\newblock In {\em Proceedings of the 31st Annual ACM/IEEE Symposium on Logic in
  Computer Science}, LICS '16, pages 565--574, New York, NY, USA, 2016. ACM.
\newblock \href {https://doi.org/10.1145/2933575.2934545}
  {\path{doi:10.1145/2933575.2934545}}.

\bibitem[Jac23]{Jack}
Tom Jack.
\newblock {$\pi_4\mathbbm{S}^3\not\cong 1$ and another Brunerie number in
  CCHM}.
\newblock Extended abstract at \emph{The Second International Conference on
  Homotopy Type Theory (HoTT 2023)}, 2023.
\newblock URL:
  \url{https://hott.github.io/HoTT-2023/abstracts/HoTT-2023_abstract_21.pdf}.

\bibitem[Jam55]{James55}
I.~M. James.
\newblock Reduced product spaces.
\newblock {\em Annals of Mathematics}, 62(1):170 -- 197, 1955.

\bibitem[Kan22a]{JamesFormalization}
Rongji Kang.
\newblock {Formalisation of the James Construction}, 2022.
\newblock URL:
  \url{https://github.com/agda/cubical/tree/master/Cubical/HITs/James}.

\bibitem[Kan22b]{BMFormalization}
Rongji Kang.
\newblock {Formalisation of the James Construction}, 2022.
\newblock URL:
  \url{https://github.com/agda/cubical/tree/master/Cubical/HITs/James}.

\bibitem[Lic14]{LicataBlog14}
Daniel~R. Licata.
\newblock Another proof that univalence implies function extensionality, 2014.
\newblock Blog post at
  \url{https://homotopytypetheory.org/2014/02/17/another-proof-that-univalence-implies-function-extensionality/}.

\bibitem[Lju20]{LjungstromMsc20}
Axel Ljungstr{\" o}m.
\newblock {Computing Cohomology in Cubical Agda}.
\newblock Master's thesis, Stockholm University, 2020.

\bibitem[Lju22]{@Axel-HoTTEST}
Axel Ljungström.
\newblock {The Brunerie Number Is -2}, 2022.
\newblock Blog post at
  \url{https://homotopytypetheory.org/2022/06/09/the-brunerie-number-is-2/}.

\bibitem[Lju24]{axel-smash}
Axel Ljungström.
\newblock {Symmetric Monoidal Smash Products in Homotopy Type Theory}, 2024.
\newblock \href {https://arxiv.org/abs/2402.03523} {\path{arXiv:2402.03523}}.

\bibitem[LLM23]{LLM23}
Thomas Lamiaux, Axel Ljungstr\"{o}m, and Anders M\"{o}rtberg.
\newblock Computing cohomology rings in cubical agda.
\newblock In {\em Proceedings of the 12th ACM SIGPLAN International Conference
  on Certified Programs and Proofs}, CPP 2023, page 239–252, New York, NY,
  USA, 2023. Association for Computing Machinery.
\newblock \href {https://doi.org/10.1145/3573105.3575677}
  {\path{doi:10.1145/3573105.3575677}}.

\bibitem[LM23]{LICS23}
Axel Ljungstr{\"{o}}m and Anders M{\"{o}}rtberg.
\newblock {Formalizing {\(\pi\)}4(S\({}^{\mbox{3}}\)) {\(\cong\)}Z/2Z and
  Computing a Brunerie Number in Cubical Agda}.
\newblock In {\em {LICS}}, pages 1--13, 2023.
\newblock \href {https://doi.org/10.1109/LICS56636.2023.10175833}
  {\path{doi:10.1109/LICS56636.2023.10175833}}.

\bibitem[LM24]{LM24}
Axel Ljungström and Anders Mörtberg.
\newblock {Computational Synthetic Cohomology Theory in Homotopy Type Theory},
  2024.
\newblock \href {https://arxiv.org/abs/2401.16336} {\path{arXiv:2401.16336}}.

\bibitem[LS13]{LicataShulman13}
Daniel~R. Licata and Michael Shulman.
\newblock {Calculating the Fundamental Group of the Circle in Homotopy Type
  Theory}.
\newblock In {\em Proceedings of the 2013 28th Annual ACM/IEEE Symposium on
  Logic in Computer Science}, LICS '13, pages 223--232, Washington, DC, USA,
  2013. IEEE Computer Society.
\newblock \href {https://doi.org/10.1109/LICS.2013.28}
  {\path{doi:10.1109/LICS.2013.28}}.

\bibitem[LS20]{LumsdaineShulman19}
Peter~LeFanu Lumsdaine and Michael Shulman.
\newblock Semantics of higher inductive types.
\newblock {\em Mathematical Proceedings of the Cambridge Philosophical
  Society}, 169(1):159--208, 2020.
\newblock \href {https://doi.org/10.1017/S030500411900015X}
  {\path{doi:10.1017/S030500411900015X}}.

\bibitem[LW24]{SteenrodSquares}
Axel Ljungström and David Wärn.
\newblock {The Steenrod Squares in HoTT Revisited}.
\newblock Extended abstract at \emph{Workshop on Homotopy Type Theory /
  Univalent Foundations (HoTT/UF24}, 2024.
\newblock URL:
  \url{https://hott-uf.github.io/2024/abstracts/HoTTUF_2024_paper_8.pdf}.

\bibitem[ML75]{MartinLof75itt}
Per Martin-L{\"o}f.
\newblock {An Intuitionistic Theory of Types: Predicative Part}.
\newblock In H.~E. Rose and J.~C. Shepherdson, editors, {\em Logic Colloquium
  '73}, volume~80 of {\em Studies in Logic and the Foundations of Mathematics},
  pages 73--118. North-Holland, 1975.
\newblock \href {https://doi.org/10.1016/S0049-237X(08)71945-1}
  {\path{doi:10.1016/S0049-237X(08)71945-1}}.

\bibitem[ML84]{MartinLof84bibliopolis}
Per Martin-L{\"o}f.
\newblock {\em Intuitionistic type theory}, volume~1 of {\em Studies in Proof
  Theory}.
\newblock Bibliopolis, 1984.

\bibitem[MP20]{cubicalsynthetic}
Anders M\"{o}rtberg and Lo\"{\i}c Pujet.
\newblock {Cubical Synthetic Homotopy Theory}.
\newblock In {\em Proceedings of the 9th ACM SIGPLAN International Conference
  on Certified Programs and Proofs}, CPP 2020, pages 158--171, New York, NY,
  USA, 2020. Association for Computing Machinery.
\newblock \href {https://doi.org/10.1145/3372885.3373825}
  {\path{doi:10.1145/3372885.3373825}}.

\bibitem[{Red}a]{cooltt}
{RedPRL Development Team}.
\newblock \texttt{\textcolor[rgb]{.012,.27,.46}{cool}tt}.
\newblock \url{https://www.github.com/RedPRL/cooltt}.

\bibitem[{Red}b]{redtt}
{RedPRL Development Team}.
\newblock \texttt{\textcolor[rgb]{.91,.31,.27}{red}tt}.
\newblock \url{https://www.github.com/RedPRL/redtt}.

\bibitem[Shu19]{Shulman19}
Michael Shulman.
\newblock All {$(\infty,1)$}-toposes have strict univalent universes, April
  2019.
\newblock Preprint.
\newblock URL: \url{https://arxiv.org/abs/1904.07004}, \href
  {https://arxiv.org/abs/1904.07004} {\path{arXiv:1904.07004}}.

\bibitem[SK22]{Syllepsis}
Kristina Sojakova and G.~A. Kavvos.
\newblock {Syllepsis in Homotopy Type Theory}.
\newblock In {\em Proceedings of the 37th Annual ACM/IEEE Symposium on Logic in
  Computer Science}, LICS '22, New York, NY, USA, 2022. Association for
  Computing Machinery.
\newblock \href {https://doi.org/10.1145/3531130.3533347}
  {\path{doi:10.1145/3531130.3533347}}.

\bibitem[{Uni}13]{HoTT13}
The {Univalent Foundations Program}.
\newblock {\em Homotopy Type Theory: Univalent Foundations of Mathematics}.
\newblock Self-published, Institute for Advanced Study, 2013.
\newblock URL: \url{https://homotopytypetheory.org/book/}.

\bibitem[VMA21]{VezzosiMortbergAbel19}
Andrea Vezzosi, Anders M{\" o}rtberg, and Andreas Abel.
\newblock {Cubical Agda: A Dependently Typed Programming Language with
  Univalence and Higher Inductive Types}.
\newblock {\em Journal of Functional Programming}, 31:e8, 2021.
\newblock \href {https://doi.org/10.1017/S0956796821000034}
  {\path{doi:10.1017/S0956796821000034}}.

\bibitem[Voe10a]{Voevodsky10cmu}
Vladimir Voevodsky.
\newblock The equivalence axiom and univalent models of type theory, February
  2010.
\newblock Notes from a talk at Carnegie Mellon University.
\newblock URL: \url{http://www.math.ias.edu/vladimir/files/CMU_talk.pdf}.

\bibitem[Voe10b]{Voevodsky10bonn}
Vladimir Voevodsky.
\newblock Univalent foundations, September 2010.
\newblock Notes from a talk in Bonn.
\newblock URL:
  \url{https://www.math.ias.edu/vladimir/sites/math.ias.edu.vladimir/files/Bonn_talk.pdf}.

\bibitem[W{\"a}r23]{david23}
David W{\"a}rn.
\newblock Eilenberg–maclane spaces and stabilisation in homotopy type theory.
\newblock {\em Journal of Homotopy and Related Structures}, 18(2):357--368, Sep
  2023.
\newblock \href {https://doi.org/10.1007/s40062-023-00330-5}
  {\path{doi:10.1007/s40062-023-00330-5}}.

\end{thebibliography}
